\documentclass[11pt]{article}

\usepackage{latexsym}
\usepackage{amsmath}
\usepackage{amssymb}
\usepackage{amsfonts}
\usepackage{amsthm}
\usepackage{amscd}
\usepackage{epsfig}
\usepackage{verbatim}
\usepackage{fancybox}
\usepackage{moreverb}
\usepackage{graphicx}
\usepackage{psfrag}

\newcommand{\R}{{\mathbb R}}

\newcommand{\D}{{\partial}}

\newtheorem{theorem}{Theorem}
\newtheorem{lemma}[theorem]{Lemma}

\newtheorem{corollary}[theorem]{Corollary}

\def\commutatif{\ar@{}[rd]|{\circlearrowleft}}

\title{Two-dimensional Finite Larmor Radius approximation in canonical gyrokinetic coordinates}

\author{E. Fr\'enod, A. Mouton}

\date{\today}

\begin{document}

\maketitle

\abstract{In this paper, we present some new results about the approximation of the Vlasov-Poisson system with a strong external magnetic field by the 2D finite Larmor radius model. The proofs within the present work are built by using two-scale convergence tools, and can be viewed as a new slant on previous works of Fr\'enod \& Sonnendr\"ucker and Bostan on the 2D finite Larmor Radius model. In a first part, we recall the physical and mathematical contexts. We also recall two main results from previous papers of Fr\'enod  \& Sonnendr\"ucker and Bostan. Then, we introduce a  set of variables which are so-called \textit{canonical gyrokinetic coordinates}, and we write the Vlasov equation in these new variables. Then, we establish some two-scale convergence and weak-* convergence results.}

\section{Introduction}

\indent Nowadays, domestic energy production by using magnetic confinement fusion (MCF) techniques is a huge technological and human challenge, as it is illustrated by the international scientific collaboration around ITER which is under construction in Cadarache (France). 
Since magnetic confinement, needed to reach nuclear fusion reaction, is a very complex physical phenomenon, the mathematical models which are linked with this plasma physics subject need to be rigorously studied from theoretical and numerical points of view. Such a work programme based on rigorous mathematical studies and high precision numerical simulations can bring some additional informations about the behavior of the studied plasma before the launch of real experiments. \\
\indent The present paper can be viewed as a part of the recent work programme about the mathematical justification of the mathematical models which are used for numerical simulations of MCF experiments. Indeed, the first tokamak plasma models have been proposed by Littlejohn, Lee \textit{et al.}, Dubin \textit{et al.} or Brizard \textit{et al.} (see \cite{Littlejohn}, \cite{Lee, Lee_2}, \cite{Dubin}, \cite{Brizard_PhD, Hahm-Brizard}) nevertheless most of these models were established by using formal assumptions. For ten years, many mathematicians have been working on mathematical justification of these models, especially the gyrokinetic approaches like guiding-center approximations and finite Larmor radius approximations: many results in this research field are due to Fr\'enod \& Sonnendr\"ucker \textit{et al.} \cite{Two-scale_expansion, Homogenization_VP, Long_time_behavior, Finite_Larmor_radius}, Golse \& Saint-Raymond \cite{Golse_1, Golse_2}, Bostan \cite{Bostan_2007} or, more recently, Han-Kwan \cite{Han-Kwan}. These mathematical results mostly rely on two-scale convergence theory (see Allaire \cite{Allaire}, Nguetseng \cite{NGuetseng}) or compactness arguments. \\
\indent In this paper, we are focused on the 2D finite Larmor radius model and its mathematical justification: more precisely, the goal is to make a synthesis of previous mathematical proofs of the convergence of $(f_{\epsilon},\tilde{\mathbf{E}}_{\epsilon})$, where 
\begin{equation}
f_{\epsilon}(\mathbf{x},k,\alpha,t) = \tilde{f}_{\epsilon}\big(x_{1}-\sqrt{2k}\,\sin\alpha, x_{2}+\sqrt{2k}\,\cos\alpha,\sqrt{2k}\,\cos\alpha, \sqrt{2k}\,\sin\alpha, t\big) \, ,
\end{equation}
and where $(\tilde{f}_{\epsilon},\tilde{\mathbf{E}}_{\epsilon})$ is the solution of the following 2D Vlasov-Poisson system
\begin{equation} \label{VP_intro}
\left\{
\begin{split}
&\D_{t} \tilde{f}_{\epsilon} + \cfrac{1}{\epsilon} \, \tilde{\mathbf{v}} \cdot \nabla_{\tilde{\mathbf{x}}} \tilde{f}_{\epsilon} + \Big(  \tilde{\mathbf{E}}_{\epsilon} + \cfrac{1}{\epsilon} \, \left(
\begin{array}{c}
\tilde{v}_{2} \\ -\tilde{v}_{1}
\end{array}
\right) \Big) \cdot \nabla_{\tilde{\mathbf{v}}} \tilde{f}_{\epsilon} = 0 \, , \\
& \tilde{f}_{\epsilon}(\tilde{\mathbf{x}},\tilde{\mathbf{v}},0) = \tilde{f}^{0}(\tilde{\mathbf{x}},\tilde{\mathbf{v}}) \, , \\
& -\nabla_{\tilde{\mathbf{x}}} \tilde{\phi}_{\epsilon} = \tilde{\mathbf{E}}_{\epsilon} \, , \qquad -\Delta_{\tilde{\mathbf{x}}} \tilde{\phi}_{\epsilon} = \int_{\R^{2}} \tilde{f}_{\epsilon} \, d\tilde{\mathbf{v}} - \tilde{n}_{e} \, ,
\end{split}
\right.
\end{equation}
towards the couple $(f,\tilde{\mathbf{E}})$ which is the solution of the 2D finite Larmor radius model given in \cite{Bostan_2007} (see also Theorem \ref{theoreme_Bostan} below). The main results on this model are due to Sonnendr\"ucker, Fr\'enod and Bostan, and indicate that $(f_{\epsilon},\tilde{\mathbf{E}}_{\epsilon})$ somehow weak-* converges to the solution of the finite Larmor radius model. However the proofs within these articles are based on various assumptions and use various tools. \\ 

\indent The first part of the present paper is devoted to a state-of-the-art about the two-dimensional finite Larmor radius approximation. Firstly, we recall the procedure which allows us to obtain the dimensionless model (\ref{VP_intro}) from the complete Vlasov-Poisson model by considering specific assumptions. Then we recall the two-scale convergence theorem of Fr\'enod \& Sonnendr\"ucker \cite{Finite_Larmor_radius} on the one hand, and the weak-* convergence theorem of Bostan \cite{Bostan_2007} on the other hand. \\
\indent In a second part, we introduce a set of variables which are so-called \textit{canonical gyrokinetic coordinates} and we reformulate the Vlasov-Poisson system (\ref{VP_intro}) in these new variables. Then, we establish a two-scale convergence theorem which only relies on Fr\'enod \& Sonnendr\"ucker's assumptions. Finally, we deduce directly from this result Bostan's Finite Larmor radius model (see \cite{Bostan_2007}) by adding a somehow non-physical assumption on the electric field which consists in considering a strong convergence of the sequence $(\tilde{\mathbf{E}}_{\epsilon})_{\epsilon\,>\,0}$.


\section{State-of-the-art}

\subsection{Scaling of the Vlasov-Poisson model}

This paragraph is devoted to the scaling of the following Vlasov-Poisson model:
\begin{equation} \label{VP-3D}
\left\{
\begin{array}{l}
\displaystyle \D_{t} \tilde{f} + \tilde{\mathbf{v}} \cdot \nabla_{\tilde{\mathbf{x}}} \tilde{f} + \frac{e}{m_{i}} \big( \tilde{\mathbf{E}} + \tilde{\mathbf{v}} \times \tilde{\mathbf{B}} \big) \cdot \nabla_{\tilde{\mathbf{v}}} \tilde{f} = 0 \, , \\
\displaystyle \tilde{f}(\tilde{\mathbf{x}},\tilde{\mathbf{v}},0) = \tilde{f}^{0}(\tilde{\mathbf{x}},\tilde{\mathbf{v}}) \, , \\
\displaystyle -\nabla_{\tilde{\mathbf{x}}} \tilde{\phi} = \tilde{\mathbf{E}} \, , \qquad -\Delta_{\tilde{\mathbf{x}}}\tilde{\phi} = \frac{e}{\varepsilon_{0}} \int_{\R_{\tilde{\mathbf{v}}}^{3}} \tilde{f}\, d\tilde{\mathbf{v}} - \frac{e}{\varepsilon_{0}} 
\, \tilde{n}_{e} \, ,
\end{array}
\right.
\end{equation}
where $\tilde{\mathbf{x}} = (\tilde{x}_{1},\tilde{x}_{2},\tilde{x}_{3}) \in \R_{\tilde{\mathbf{x}}}^{3}$ is the position variable, $\tilde{\mathbf{v}} = (\tilde{v}_{1},\tilde{v}_{2},\tilde{v}_{3}) \in \R_{\tilde{\mathbf{v}}}^{3}$ is the velocity variable, $t \in \R_{+}$ is the time variable, $\tilde{f} = \tilde{f}(\tilde{\mathbf{x}},\tilde{\mathbf{v}},t)$ is the ion distribution function, $\tilde{n}_{e}$ is the electron density, $\tilde{\mathbf{E}} = \tilde{\mathbf{E}}(\tilde{\mathbf{x}},t)$ is the self-consistent electric field generated by the ions and the electrons, $\tilde{\mathbf{B}} = \tilde{\mathbf{B}}(\tilde{\mathbf{x}},t)$ is the magnetic field which is applied on the considered plasma, $\tilde{\phi} = \tilde{\phi}(\tilde{\mathbf{x}},t)$ is the electric potential linked with $\tilde{\mathbf{E}}$, $e$ is the elementary charge and $m_{i}$ is the elementary mass of an ion. \\

\indent In this model, the external magnetic field $\tilde{\mathbf{B}}$ is assumed to be uniform and carried by the unit vector $\mathbf{e}_{3}$. We also assume that the electron density $\tilde{n}_{e}$ is given for any $(\tilde{\mathbf{x}},t) \in \R_{\tilde{\mathbf{x}}}^{3} \times \R_{+}$. Following the same approach as in Bostan \cite{Bostan_2007}, Fr\'enod \textit{et al.} \cite{Two-scale_expansion, Finite_Larmor_radius}, Golse \textit{et al.} \cite{Golse_1, Golse_2} and Han-Kwan \cite{Han-Kwan}, we add the following assumptions:
\begin{itemize}
\item[(i)] The magnetic field is supposed to be strong,
\item[(ii)] The finite Larmor radius effects are taken into account,
\item[(iii)] The ion gyroperiod is supposed to be small.
\end{itemize}

We define the dimensionless variables and unknowns $\tilde{\mathbf{x}}'=(\tilde{x}_{1}',\tilde{x}_{2}',\tilde{x}_{3}')$, $\tilde{\mathbf{v}}'=(\tilde{v}_{1}',\tilde{v}_{2}',\tilde{v}_{3}')$, $t'$, $\tilde{f}'$, $\tilde{\mathbf{E}}'$ and $\tilde{\phi}'$ by
\begin{equation} \label{scaling_def}
\begin{array}{c}
\tilde{x}_{1} = \overline{L_{\perp}} \, \tilde{x}_{1}' \, , \qquad \tilde{x}_{2} = \overline{L_{\perp}} \, \tilde{x}_{2}' \, , \qquad \tilde{x}_{3} = \overline{L_{||}} \, \tilde{x}_{3}' \, , \qquad t = \overline{t} \, t' \, , \qquad \tilde{\mathbf{v}} = \overline{v} \, \tilde{\mathbf{v}}' \, , \\ \\
\tilde{f}(\overline{L_{\perp}}\,\tilde{x}_{1}',\overline{L_{\perp}}\,\tilde{x}_{2}', \overline{L_{||}}\,\tilde{x}_{3}', \overline{v}\,\tilde{v}_{1}', \overline{v}\,\tilde{v}_{2}', \overline{v}\,\tilde{v}_{3}',\overline{t}\,t') 
=  \overline{f} \, \tilde{f}'(\tilde{x}_{1}',\tilde{x}_{2}',\tilde{x}_{3}',\tilde{v}_{1}',\tilde{v}_{2}',\tilde{v}_{3}',t') \, , \\ \\
\tilde{\mathbf{E}}(\overline{L_{\perp}}\,\tilde{x}_{1}',\overline{L_{\perp}}\,\tilde{x}_{2}', \overline{L_{||}}\,\tilde{x}_{3}', \overline{t}\,t') = \overline{E}\, \tilde{\mathbf{E}}'(\tilde{x}_{1}',\tilde{x}_{2}',\tilde{x}_{3}',t') \, , \\ \\
\tilde{\phi}(\overline{L_{\perp}}\,\tilde{x}_{1}',\overline{L_{\perp}}\,\tilde{x}_{2}', \overline{L_{||}}\,\tilde{x}_{3}', \overline{t}\,t') = \overline{\phi}\, \tilde{\phi}'(\tilde{x}_{1}',\tilde{x}_{2}',\tilde{x}_{3}',t') \, .
\end{array}
\end{equation}
In these definitions, $\overline{L_{\perp}}$ is the characteristic length in the direction perpendicular to the magnetic field, $\overline{L_{||}}$ is the characteristic length in the direction of the magnetic field, $\overline{v}$ is the characteristic velocity and $\overline{t}$ is the characteristic time. 
We also rescale the electron density as follows:
\begin{equation}
\tilde{n}_{e}(\overline{L_{\perp}}\,\tilde{x}_{1}',\overline{L_{\perp}}\,\tilde{x}_{2}', \overline{L_{||}}\,\tilde{x}_{3}', \overline{t}\,t') = \overline{n} \, \tilde{n}_{e}'(\tilde{x}_{1}',\tilde{x}_{2}',\tilde{x}_{3}',t') \, .
\end{equation}
\indent Following the assumptions on the magnetic field $\tilde{\mathbf{B}}$, we set $\overline{B}$ as being such that
\begin{equation}
\tilde{\mathbf{B}} = \overline{B} \, \mathbf{e}_{3} \, .
\end{equation}
Then, we set  $ \overline{L_{||}}$ as the size of the physical device in game.
We also link $\overline{f}$, $\overline{E}$ and $\overline{\phi}$ with the characteristic Debye length $\overline{\lambda_{D}}$ by
\begin{equation}
\overline{f} = \cfrac{\overline{n}}{\overline{v}^{3}} \, , \qquad \overline{E} = \cfrac{\overline{\lambda_{D}} \, e\, \overline{n}}{\varepsilon_{0}} \, , \qquad \overline{\phi} = \cfrac{\overline{\lambda_{D}}^{2} \, e \, \overline{n}}{\varepsilon_{0}} \, ,
\end{equation}
and we take $\overline{\lambda_{D}}$ as the characteristic length in the direction perpendicular to the magnetic field, \textit{i.e.}
\begin{equation}
\overline{L_{\perp}} = \overline{\lambda_{D}} \, .
\end{equation}
Since we want to take into account the smallness of the gyroperiod and the finite Larmor radius effects, we  define the characteristic gyrofrequency $\overline{\omega_{i}}$  and  the characteristic Larmor radius $\overline{r_{L}}$ as
\begin{equation}
\overline{\omega_{i}} = \cfrac{e\,\overline{B}}{m_{i}}, \qquad 
\overline{r_{L}} = \cfrac{\overline{v}}{\overline{\omega_{i}}} \, .
\end{equation}
\indent With these notations, the Vlasov-Poisson system is rescaled as follows:
\begin{equation} \label{VP-3D_rescaled_1}
\left\{
\begin{array}{l}
\displaystyle \D_{t'}\tilde{f}' + \overline{t}\,\overline{v} \, \left(
\begin{array}{c}
\frac{\overline{r_{L}}}{\overline{\lambda_{D}}} \, \tilde{v}_{1}' \\ \frac{\overline{r_{L}}}{\overline{\lambda_{D}}} \, \tilde{v}_{2}' \\ \frac{\overline{r_{L}}}{\overline{L_{||}}} \, \tilde{v}_{3}'
\end{array}
\right) \cdot \nabla_{\tilde{\mathbf{x}}'}\tilde{f}' + \Big( \cfrac{\overline{t} \, \overline{\lambda_{D}} \, e^{2}\, \overline{n}}{\varepsilon_{0}\,m_{i}\,\overline{v}} \, \tilde{\mathbf{E}}' + \overline{t} \, \overline{\omega_{i}} \, \tilde{\mathbf{v}}' \times \mathbf{e}_{3} \Big) \cdot \nabla_{\tilde{\mathbf{v}}'} \tilde{f}' = 0 \, , \\ \\
\tilde{f}'(\tilde{\mathbf{x}}',\tilde{\mathbf{v}}',0) = \tilde{f^{0}}'(\tilde{\mathbf{x}}', \tilde{\mathbf{v}}') \, , \\ \\
\tilde{\mathbf{E}}' = \left(
\begin{array}{c}
- \D_{\tilde{x}_{1}'} \tilde{\phi}' \\ -\D_{\tilde{x}_{2}'} \tilde{\phi}' \\ -\frac{\overline{\lambda_{D}}}{\overline{L_{||}}}\, \D_{\tilde{x}_{3}'}\tilde{\phi}'
\end{array}
\right) \, , \displaystyle \qquad -\Delta_{(\tilde{x}_{1}',\tilde{x}_{2}')}\tilde{\phi}'- \cfrac{\overline{\lambda_{D}}^{2}}{\overline{L_{||}}^{2}} \, \D_{\tilde{x}_{3}'}^{2}\tilde{\phi}' = \int_{\R^{3}} \tilde{f}' \, d\tilde{\mathbf{v}}' - \tilde{n}_{e}' \, .
\end{array}
\right.
\end{equation}

\indent Taking into account the finite Larmor radius effects consists in considering a regime in which the 
Larmor radius is of the order of the Debye length. This implies
\begin{equation}
\cfrac{\overline{r_{L}}}{\overline{\lambda_{D}}} = 1 \, .
\end{equation}
Since the magnetic field is assumed to be strong, the Larmor radius is small when compared with the size of the physical domain. Then it is natural to take
\begin{equation}
\cfrac{\overline{r_{L}}}{\overline{L_{||}}} = \epsilon \, ,
\end{equation}
where $\epsilon > 0$ is small. \\
Assumption (iii) can be translated in terms of characteristic scales by
\begin{equation}
\overline{t} \, \overline{\omega_{i}} = \cfrac{1}{\epsilon} \, .
\end{equation}
Assumption (i) means that the magnetic force is much stronger than the electric force, so we consider
\begin{equation}
\cfrac{\overline{E} \, e}{\overline{v} \, m_{i}\, \overline{\omega}_{i}} = \epsilon \, .
\end{equation}
Then, removing the primes and adding $\epsilon$ in subscript, the rescaled Vlasov-Poisson model writes
\begin{equation} \label{VP-3D_rescaled}
\left\{
\begin{split}
&\D_{t} \tilde{f}_{\epsilon} + \cfrac{1}{\epsilon} \, \left(
\begin{array}{c}
\tilde{v}_{1} \\ \tilde{v}_{2}
\end{array}
\right) \cdot \nabla_{(\tilde{x}_{1},\tilde{x}_{2})} \tilde{f}_{\epsilon} +  \tilde{v}_{3} \, \D_{\tilde{x}_{3}}\tilde{f}_{\epsilon} + \Big( \tilde{\mathbf{E}}_{\epsilon} + \cfrac{1}{\epsilon} \, \tilde{\mathbf{v}} \times \mathbf{e}_{3} \Big) \cdot \nabla_{\tilde{\mathbf{v}}} \tilde{f}_{\epsilon} = 0 \, , \\
& \tilde{f}_{\epsilon}(\tilde{\mathbf{x}},\tilde{\mathbf{v}},0) = \tilde{f}_{\epsilon}^{0}(\tilde{\mathbf{x}},\tilde{\mathbf{v}}) \, , \\
&\left(
\begin{array}{c}
-\nabla_{(\tilde{x}_{1},\tilde{x}_{2})} \tilde{\phi}_{\epsilon} \\ - \epsilon \, \D_{\tilde{x}_{3}}\tilde{\phi}_{\epsilon}
\end{array}
\right) = \tilde{\mathbf{E}}_{\epsilon} \, , \qquad -\Delta_{(\tilde{x}_{1},\tilde{x}_{2})} \tilde{\phi}_{\epsilon} - \epsilon^{2} \, \D_{\tilde{x}_{3}}^{2}\tilde{\phi}_{\epsilon} = \int_{\R^{3}} \tilde{f}_{\epsilon} \, d\tilde{\mathbf{v}} - \tilde{n}_{e} \, ,
\end{split}
\right.
\end{equation}
which is the model studied in previous works of Fr\'enod \& Sonnendr\"ucker \cite{Finite_Larmor_radius}, Golse \& Saint-Raymond \cite{Golse_1, Golse_2}, and Bostan \cite{Bostan_2007}.

\subsection{Previous results}

In this paragraph, we recall two main results about the asymptotic behavior of the sequences $(\tilde{f}_{\epsilon})_{\epsilon \, > \, 0}$ and $(\tilde{\mathbf{E}}_{\epsilon})_{\epsilon \, > \, 0}$ when $\epsilon$ goes to 0. The first one is based on the use of two-scale convergence and homogenization techniques developed by Allaire \cite{Allaire} and Nguetseng \cite{NGuetseng}, and was established by Fr\'enod and Sonnendr\"ucker in \cite{Finite_Larmor_radius}. The second one relies on compactness arguments and was proved by Bostan in \cite{Bostan_2007}. After recalling these two results, we discuss the main differences between them. These differences are the source of the motivation of the present paper. \\

\indent In order to simplify, we consider that the whole model (\ref{VP-3D_rescaled}) does not depend on $\tilde{x}_{3}$ nor $\tilde{v}_{3}$, and we assume that $\tilde{f}_{\epsilon}^{0} = \tilde{f}^{0}$ for all $\epsilon$. Then, it is reduced to a singularly perturbed 2D Vlasov-Poisson model of the form
\begin{equation} \label{VP-2D_rescaled}
\left\{
\begin{split}
&\D_{t} \tilde{f}_{\epsilon} + \cfrac{1}{\epsilon} \, \tilde{\mathbf{v}} \cdot \nabla_{\tilde{\mathbf{x}}} \tilde{f}_{\epsilon} + \Big(  \tilde{\mathbf{E}}_{\epsilon} + \cfrac{1}{\epsilon} \, \left(
\begin{array}{c}
\tilde{v}_{2} \\ -\tilde{v}_{1}
\end{array}
\right) \Big) \cdot \nabla_{\tilde{\mathbf{v}}} \tilde{f}_{\epsilon} = 0 \, , \\
& \tilde{f}_{\epsilon}(\tilde{\mathbf{x}},\tilde{\mathbf{v}},0) = \tilde{f}^{0}(\tilde{\mathbf{x}},\tilde{\mathbf{v}}) \, , \\
& -\nabla_{\tilde{\mathbf{x}}} \tilde{\phi}_{\epsilon} = \tilde{\mathbf{E}}_{\epsilon} \, , \qquad -\Delta_{\tilde{\mathbf{x}}} \tilde{\phi}_{\epsilon} = \int_{\R^{2}} \tilde{f}_{\epsilon} \, d\tilde{\mathbf{v}} - \tilde{n}_{e} \, ,
\end{split}
\right.
\end{equation}
where $\tilde{\mathbf{x}} = (\tilde{x}_{1},\tilde{x}_{2}) \in \R^{2}$ and $\tilde{\mathbf{v}} = (\tilde{v}_{1},\tilde{v}_{2}) \in \R^{2}$. 

~

The following theorem can be attributed to Fr\'enod and Sonnendr\"ucker \cite{Finite_Larmor_radius}
(see Theorem 1.5) even if the setting of this paper is a charged particle beam Vlasov-Poisson model not involving any electron density. Nonetheless,  the proof of \cite{Finite_Larmor_radius} works again in the setting of  model (\ref{VP-2D_rescaled})
which involves an electron density.

\begin{theorem}[Fr\'enod \& Sonnendr\"ucker \cite{Finite_Larmor_radius}] \label{CV_cart}
We assume that, for a fixed $p \geq 2$, $\tilde{f}^{0}$ is in $L^{1}(\R^{4}) \cap L^{p}(\R^{4})$, is positive everywhere and such that 
\begin{equation}
\displaystyle \int_{\R^{4}} |\tilde{\mathbf{v}}|^{2}\,\tilde{f}^{0}(\tilde{\mathbf{x}},\tilde{\mathbf{v}}) \, d\tilde{\mathbf{x}}\,d\tilde{\mathbf{v}} < +\infty \, .
\end{equation}
We also assume that $\tilde{n}_{e}$ does not depend on $t$, is in $L^{1}(\R^{2})\cap L^{3/2}(\R^{2})$ and satisfies 
\begin{equation}
\displaystyle \int_{\R^{4}} \tilde{f}^{0}(\tilde{\mathbf{x}},\tilde{\mathbf{v}}) \, d\tilde{\mathbf{x}} \, d\tilde{\mathbf{v}} = \int_{\R^{2}}\tilde{n}_{e}(\tilde{\mathbf{x}}) \, d\tilde{\mathbf{x}} \, .
\end{equation}
Then, the sequence $(\tilde{f}_{\epsilon},\tilde{\mathbf{E}}_{\epsilon})_{\epsilon\,>\,0}$ is bounded $L^{\infty}\big(0,T;L^{p}(\R^{4})\big) \times \big(L^{\infty}\big(0,T;W^{1,3/2}(\R^{2})\big)\big)^{2}$ independently of $\epsilon$. Furthermore, by extracting some subsequences, 
\begin{equation}
\begin{array}{ccll}
\tilde{f}_{\epsilon} &\longrightarrow & \tilde{F} = \tilde{F}(\tilde{\mathbf{x}},\tilde{\mathbf{v}},\tau,t) & \textit{two-scale in $L^{\infty}\big(0,T;L_{\#}^{\infty}\big(0,2\pi; L^{p}(\R^{4})\big)\big)$,} \\
\tilde{\mathbf{E}}_{\epsilon} &\longrightarrow & \tilde{\mathcal{E}}=\tilde{\mathcal{E}}(\tilde{\mathbf{x}},\tau,t) & \textit{two-scale in $\Big(L^{\infty}\big(0,T;L_{\#}^{\infty}\big(0,2\pi;W^{1,3/2}(\R^{2})\big)\big)\Big)^{2}$.}
\end{array}
\end{equation}
Moreover, $\tilde{F}$ is linked with $\tilde{G} = \tilde{G}(\tilde{\mathbf{y}},\tilde{\mathbf{u}},t) \in L^{\infty}\big(0,T;L^{p}(\R^{4})\big)$ by the relation
\begin{equation}
\tilde{F}(\tilde{\mathbf{x}},\tilde{\mathbf{v}},\tau,t) = \tilde{G}\big(\tilde{\mathbf{x}}+\mathcal{R}(-\tau)\,\tilde{\mathbf{v}},R(-\tau)\,\tilde{\mathbf{v}},t\big) \, ,
\end{equation}
and $(\tilde{G},\tilde{\mathcal{E}})$ is the solution of
\begin{equation} \label{H-cart}
\left\{
\begin{array}{l}
\displaystyle \D_{t} \tilde{G}(\tilde{\mathbf{y}},\tilde{\mathbf{u}},t) + \Bigg[\int_{0}^{2\pi} \mathcal{R}(-\sigma) \, \tilde{\mathcal{E}}\big(\tilde{\mathbf{y}}+\mathcal{R}(\sigma)\,\tilde{\mathbf{u}},\sigma,t\big) \, d\sigma \Bigg] \cdot \nabla_{\tilde{\mathbf{y}}} \tilde{G}(\tilde{\mathbf{y}},\tilde{\mathbf{u}},t) \\
\displaystyle \qquad \qquad \qquad + \Bigg[\int_{0}^{2\pi} R(-\sigma) \, \tilde{\mathcal{E}}\big(\tilde{\mathbf{y}}+\mathcal{R}(\sigma)\,\tilde{\mathbf{u}},\sigma,t\big) \, d\sigma \Bigg] \cdot \nabla_{\tilde{\mathbf{u}}} \tilde{G}(\tilde{\mathbf{y}},\tilde{\mathbf{u}},t) = 0 \, , \\
\displaystyle \tilde{G}(\tilde{\mathbf{y}},\tilde{\mathbf{u}},0) = \frac{1}{2\pi} \, \tilde{f}^{0}(\tilde{\mathbf{y}},\tilde{\mathbf{u}}) \, , \\ \\ 
\displaystyle \tilde{\mathcal{E}}(\tilde{\mathbf{x}},\tau,t) = -\nabla_{\tilde{\mathbf{x}}} \tilde{\Phi}(\tilde{\mathbf{x}},\tau,t) \, , \\ \\
\displaystyle -\Delta_{\tilde{\mathbf{x}}} \tilde{\Phi}(\tilde{\mathbf{x}},\tau,t) = \int_{\R^{2}} \tilde{G}\big(\tilde{\mathbf{x}}+\mathcal{R}(-\tau)\,\tilde{\mathbf{v}}, R(-\tau)\,\tilde{\mathbf{v}},t\big) \, d\tilde{\mathbf{v}} - \frac{1}{2\pi}\,\tilde{n}_{e}(\tilde{\mathbf{x}})\, ,
\end{array}
\right.
\end{equation}
where
\begin{equation} \label{mat_rotation}
\mathcal{R}(\tau) = \left(
\begin{array}{cc}
\sin\tau & 1-\cos\tau \\
\cos\tau -1 & \sin\tau
\end{array}
\right) \, , \quad R(\tau) = \left(
\begin{array}{cc}
\cos\tau & \sin\tau \\
-\sin\tau & \cos\tau
\end{array}
\right) \, .
\end{equation}
\end{theorem}

In this theorem, $L_{\#}^{\infty}(0,2\pi; L^{p}(\R^4))$ stands for the space of functions $\tilde{h} = \tilde{h}(\tilde{\mathbf{x}},\tilde{\mathbf{v}},\tau)$ being in $L^{\infty}(0,2\pi; L^{p}(\R^4))$ and $2\pi$-periodic with respect to $\tau$.

\indent As a consequence of this theorem, extracting some subsequences, we have
\begin{equation} \label{CV_weak_f}
\begin{array}{cccl}
\tilde{f}_{\epsilon} & \stackrel{*}{\rightharpoonup} & \tilde{f} & \textnormal{in $L^{\infty}\big(0,T;L^{p}(\R^{4})\big)$,} \\
\tilde{\mathbf{E}}_{\epsilon} & \stackrel{*}{\rightharpoonup} & \tilde{\mathbf{E}} & \textnormal{in $\big(L^{\infty}\big(0,T;W^{1,3/2}(\R^{2})\big)\big)^{2}$,}
\end{array}
\end{equation}
where
\begin{equation}
\tilde{f}(\tilde{\mathbf{x}},\tilde{\mathbf{v}},t) = \int_{0}^{2\pi} \tilde{F}(\tilde{\mathbf{x}},\tilde{\mathbf{v}},\tau,t)\,d\tau \, , \qquad  \textnormal{and} \qquad \tilde{\mathbf{E}}(\tilde{\mathbf{x}},t) = \int_{0}^{2\pi} \tilde{\mathcal{E}}(\tilde{\mathbf{x}},\tau,t)\,d\tau \, .
\end{equation}
By using the relation between $\tilde{F}$ and $\tilde{G}$, we can easily remark that $(\tilde{f},\tilde{\mathbf{E}})$ is solution of
\begin{equation}
\left\{
\begin{array}{l}
\displaystyle \D_{t}\tilde{f}(\tilde{\mathbf{x}}, \tilde{\mathbf{v}},t) \\
\qquad \displaystyle + \int_{0}^{2\pi} \Bigg[\int_{0}^{2\pi} \mathcal{R}(\tau-\sigma) \, \tilde{\mathcal{E}}\big(\tilde{\mathbf{x}}+\mathcal{R}(\sigma-\tau)\,\tilde{\mathbf{v}},\sigma,t\big) \, d\sigma \Bigg] \cdot \nabla_{\tilde{\mathbf{x}}} \tilde{F}(\tilde{\mathbf{x}},\tilde{\mathbf{v}},\tau,t) \, d\tau \\
\qquad \displaystyle + \int_{0}^{2\pi} \Bigg[\int_{0}^{2\pi} R(\tau-\sigma) \, \tilde{\mathcal{E}}\big(\tilde{\mathbf{x}}+\mathcal{R}(\sigma-\tau)\,\tilde{\mathbf{v}},\sigma,t\big) \, d\sigma \Bigg] \cdot \nabla_{\tilde{\mathbf{v}}} \tilde{F}(\tilde{\mathbf{x}},\tilde{\mathbf{v}},\tau,t) \, d\tau = 0 \, , \\
\displaystyle \tilde{f}(\tilde{\mathbf{x}},\tilde{\mathbf{v}},0) = \frac{1}{2\pi} \int_{0}^{2\pi} \tilde{f}^{0}\big(\tilde{\mathbf{x}}+\mathcal{R}(-\tau)\,\tilde{\mathbf{v}},R(-\tau)\,\tilde{\mathbf{v}}\big) \, d\tau \, , \\
-\nabla_{\tilde{\mathbf{x}}}\tilde{\phi}(\tilde{\mathbf{x}},t) = \tilde{\mathbf{E}}(\tilde{\mathbf{x}},t) \, , \qquad  -\Delta_{\tilde{\mathbf{x}}}\tilde{\phi}(\tilde{\mathbf{x}},t) = \displaystyle\int_{\R^{2}} \tilde{f}(\tilde{\mathbf{x}},\tilde{\mathbf{v}},t) \, d\tilde{\mathbf{v}} - \tilde{n}_{e}(\tilde{\mathbf{x}}) \, .
\end{array}
\right.
\end{equation}
We notice that these equations still involve $\tilde{F}$ and $\tilde{\mathcal{E}}$. In order to make these dependencies disappear, Bostan has proposed in \cite{Bostan_2007} a reformulation of the Vlasov equation in guiding-center coordinates. Before presenting it, we introduce the sequence $(\breve{f}_{\epsilon})_{\epsilon\,> \,0}$ defined by
\begin{equation}
\tilde{f}_{\epsilon}(\tilde{\mathbf{x}},\tilde{\mathbf{v}},t) = \breve{f}_{\epsilon}\Big( \tilde{\mathbf{x}} + \left(
\begin{array}{c}
\tilde{v}_{2} \\ -\tilde{v}_{1}
\end{array}
\right), \tilde{\mathbf{v}},t\Big) \, ,
\end{equation}
and, in the same spirit, we define the initial guiding-center distribution $\breve{f}^{0}$ by
\begin{equation}
\tilde{f}^{0}(\tilde{\mathbf{x}},\tilde{\mathbf{v}},t) = \breve{f}^{0}\Big( \tilde{\mathbf{x}} + \left(
\begin{array}{c}
\tilde{v}_{2} \\ -\tilde{v}_{1}
\end{array}
\right), \tilde{\mathbf{v}},t\Big) \, .
\end{equation}

\begin{theorem}[Bostan \cite{Bostan_2007}] \label{theoreme_Bostan}
We assume that $\tilde{n}_{e} = 1$, and that $\tilde{f}^{0}$ is $2\pi$-periodic in $\tilde{x}_{1}$ and $\tilde{x}_{2}$, is positive everywhere and satisfies
\begin{equation} \label{f0_somme1}
\begin{split}
&\int_{\R^{2}}\int_{0}^{2\pi}\int_{0}^{2\pi} \tilde{f}^{0}(\tilde{\mathbf{x}},\tilde{\mathbf{v}}) \, d\tilde{\mathbf{x}} \, d\tilde{\mathbf{v}} = 1 \, , \quad \int_{\R^{2}} \int_{0}^{2\pi}\int_{0}^{2\pi} |\tilde{\mathbf{v}}|^{2}\,\tilde{f}^{0}(\tilde{\mathbf{x}},\tilde{\mathbf{v}}) \, d\tilde{\mathbf{x}} \, d\tilde{\mathbf{v}} < +\infty \, .
\end{split}
\end{equation}
We also assume that there exists $\tilde{F}^{0} \in L^{\infty}(\R_{+}) \cap L^{1}(\R_{+}; r\, dr)$ such that
\begin{equation} \label{F0_decroissante}
\forall\,(\tilde{\mathbf{x}},\tilde{\mathbf{v}}) \in [0,2\pi]^{2} \times \R^{2} \, , \qquad \tilde{f}^{0}(\tilde{\mathbf{x}},\tilde{\mathbf{v}}) \leq \tilde{F}^{0}\big(|\tilde{\mathbf{v}}|\big) \, .
\end{equation}
We also assume that $(\tilde{\mathbf{E}}_{\epsilon})_{\epsilon\,>\,0}$ admits a strong limit denoted with $\tilde{\mathbf{E}}$ in $\big(L^{2}\big(0,T;L_{\#}^{2}\big([0,2\pi]^{2}\big)\big)\big)^{2}$. Then, up to a subsequence, $\breve{f}_{\epsilon}$ weakly-* converges to a function $\breve{f} = \breve{f}(\mathbf{x},\mathbf{v},t)$ in $L^{\infty}\big([0,T)\times\R^{2}; L_{\#}^{\infty}\big([0,2\pi]^{2}\big)\big)$ verifying
\begin{equation} \label{modele_Bostan_1}
\breve{f}(\mathbf{x},\mathbf{v},t) = \frac{1}{2\pi} \, g\big(\mathbf{x},\frac{|\mathbf{v}|^{2}}{2},t\big) \, ,
\end{equation}
where $g = g(\mathbf{x},k,t)$ is the solution of
\begin{equation} \label{modele_Bostan_2}
\left\{
\begin{split}
&\D_{t}g + \langle \mathcal{E}_{2} \rangle \, \D_{x_{1}}g - \langle \mathcal{E}_{1} \rangle \, \D_{x_{2}}g = 0 \, , \\
&g(\mathbf{x},k,0) = \int_{0}^{2\pi} \breve{f}^{0}(\mathbf{x},\sqrt{2k}\,\cos\alpha,\sqrt{2k}\,\sin\alpha)\, d\alpha \, , \\
&\langle \mathcal{E}\rangle(\mathbf{x},k,t) = \frac{1}{2\pi} \int_{0}^{2\pi} \tilde{\mathbf{E}}\big(x_{1}-\sqrt{2k}\,\sin\alpha,x_{2}+\sqrt{2k}\,\cos\alpha,t\big) \, d\alpha \, , \\
& -\nabla_{\tilde{\mathbf{x}}} \tilde{\phi}(\tilde{\mathbf{x}},t) = \tilde{\mathbf{E}}(\tilde{\mathbf{x}},t) \, , \\ 
&-\Delta_{\tilde{\mathbf{x}}}\tilde{\phi}(\tilde{\mathbf{x}},t) = \frac{1}{2\pi} \int_{0}^{+\infty} \int_{0}^{2\pi} g\big(\tilde{x}_{1}+\sqrt{2k}\,\sin\alpha,\tilde{x}_{2}-\sqrt{2k}\,\cos\alpha,k,t\big) \, d\alpha\,dk - 1 \, .
\end{split}
\right.
\end{equation}
\end{theorem}

In this theorem, $L_{\#}^{2}\big([0,2\pi]^{2}\big)$ stands for the space of functions $\tilde{h}=\tilde{h}(\tilde{\mathbf{x}})$ being in $L^{2}\big([0,2\pi]^{2}\big)$ and $2\pi$-periodic with respect to $\tilde{x}_{1}$ and $\tilde{x}_{2}$.

~

\indent This last result introduces a mathematical justification of the approximation of the Vlasov-Poisson model (\ref{VP-2D_rescaled}) by the finite Larmor radius model which is exactly (\ref{modele_Bostan_1})-(\ref{modele_Bostan_2}). However, in order to prove this convergence result, Bostan considered stronger assumptions on $\tilde{f}^{0}$ and $\tilde{n}_{e}$ than  needed to get existence of the weak-* limit $(\tilde{f},\tilde{\mathbf{E}})$ from Theorem \ref{CV_cart}:  the initial distribution $\tilde{f}^{0}$ is supposed to be $2\pi$-periodic in $\tilde{x}_{1}$ and $\tilde{x}_{2}$ and $(\tilde{\mathbf{E}}_{\epsilon})_{\epsilon\,>\,0}$ is supposed to admit a strong limit in some Banach space.

\section{Convergence result in canonical gyrokinetic coordinates}

\setcounter{equation}{0}

This section is devoted to a two-scale convergence result for model (\ref{VP-2D_rescaled}) written in a new set of variables which are so-called \textit{canonical gyrokinetic coordinates}. Then, by adding a non-physical hypothesis for the electric field $\tilde{\mathbf{E}}_{\epsilon}$, we obtain straightforwardly Bostan's weak-* convergence result.


\subsection{Reformulation of Vlasov equation}

Following the ideas of Littlejohn \cite{Littlejohn}, Lee \cite{Lee, Lee_2}, and Brizard \textit{et al.} \cite{Brizard_PhD, Hahm-Brizard}, we define the variables $(x_{1},x_{2},k,\alpha) \in \R^{2} \times \R_{+} \times [0,2\pi]$ by linking them with  $(\tilde{x}_{1},\tilde{x}_{2},\tilde{v}_{1},\tilde{v}_{2}) \in \R^{4}$ by 
\begin{equation}
\left\{
\begin{array}{rcl}
\tilde{x}_{1} &=& x_{1} - \tilde{v}_{2} \, , \\
\tilde{x}_{2} &=& x_{2} + \tilde{v}_{1} \, ,
\end{array}
\right. \qquad 
\begin{array}{rcl}
\tilde{v}_{1} &=& \sqrt{2k}\,\cos\alpha \, , \\
\tilde{v}_{2} &=& \sqrt{2k}\,\sin\alpha \, .
\end{array}
\end{equation}
This set of variables is so-called \textit{canonical gyrokinetic coordinates}: indeed, if we define the characteristics $X_{1},X_{2},K,A$ linked with $x_{1},x_{2},k,\alpha$ by
\begin{equation}
\left\{
\begin{array}{rcl}
\tilde{X}_{1} &=& X_{1} - \tilde{V}_{2} \, , \\
\tilde{X}_{2} &=& X_{2} + \tilde{V}_{1} \, ,
\end{array}
\right. \qquad 
\begin{array}{rcl}
\tilde{V}_{1} &=& \sqrt{2K}\,\cos A \, , \\
\tilde{V}_{2} &=& \sqrt{2K}\,\sin A \, ,
\end{array}
\end{equation}
where $\tilde{X}_{1},\tilde{X}_{2},\tilde{V}_{1},\tilde{V}_{2}$ are the characteristics associated with the Vlasov equation (\ref{VP-2D_rescaled}.a), \textit{i.e.} satisfying
\begin{equation}
\left\{
\begin{array}{rcl}
\D_{t}\tilde{X}_{1}(t) &=& \displaystyle \frac{1}{\epsilon}\,\tilde{V}_{1}(t) \, , \\ \\
\D_{t}\tilde{X}_{2}(t) &=& \displaystyle \frac{1}{\epsilon}\,\tilde{V}_{2}(t) \, ,
\end{array}
\right. \qquad 
\begin{array}{rcl}
\D_{t}\tilde{V}_{1}(t) &=& \displaystyle \tilde{E}_{\epsilon,1}\big(\tilde{X}_{1}(t), \tilde{X}_{2}(t), t\big) + \frac{1}{\epsilon}\,\tilde{V}_{2}(t) \, , \\ \\
\D_{t}\tilde{V}_{2}(t) &=& \displaystyle \tilde{E}_{\epsilon,2}\big(\tilde{X}_{1}(t), \tilde{X}_{2}(t), t\big) - \frac{1}{\epsilon}\,\tilde{V}_{1}(t) \, ,
\end{array}
\end{equation}
we have
\begin{equation}
\left\{
\begin{array}{rcl}
\D_{t}X_{1}(t) &=& -\D_{x_{2}}H_{\epsilon}\big(X_{1}(t),X_{2}(t),K(t),A(t),t\big) \, , \\ 
\D_{t}X_{2}(t) &=& \D_{x_{1}}H_{\epsilon}\big(X_{1}(t),X_{2}(t),K(t),A(t),t\big) \, , \\
\D_{t}K(t) &=& \D_{\alpha}H_{\epsilon}\big(X_{1}(t),X_{2}(t),K(t),A(t),t\big) \, ,\\
\D_{t}A(t) &=& -\D_{k}H_{\epsilon}\big(X_{1}(t),X_{2}(t),K(t),A(t),t\big) \, ,
\end{array}
\right.
\end{equation}
where the hamiltonian function $H_{\epsilon}$ is defined by
\begin{equation}
H_{\epsilon}(x_{1},x_{2},k,\alpha,t) = \frac{k}{\epsilon} + \phi_{\epsilon}(x_{1},x_{2},k,\alpha,t) \, ,
\end{equation}
and $\phi_{\epsilon}$ is linked with $\tilde{\phi}_{\epsilon}$ by the relation
\begin{equation} \label{def_phieps_phiepstilde}
\phi_{\epsilon}(x_{1},x_{2},k,\alpha,t) = \tilde{\phi}_{\epsilon}\big(x_{1}-\sqrt{2k}\,\sin\alpha, x_{2}+\sqrt{2k}\,\cos\alpha,t\big) \, .
\end{equation}
Then it is straightforward to see that, in the  gyrokinetic canonical coordinates, the Vlasov-Poisson system (\ref{VP-2D_rescaled}) has the following shape:

\begin{equation} \label{VP2D_GC} 
\hspace{-0.2cm} \left\{
\begin{array}{l}
\displaystyle \D_{t}f_{\epsilon} + E_{\epsilon,2} \, \D_{x_{1}}f_{\epsilon} - E_{\epsilon,1}\,\D_{x_{2}}f_{\epsilon} + \sqrt{2k} \, (E_{\epsilon,1} \, \cos\alpha + E_{\epsilon,2} \, \sin\alpha) \, \D_{k}f_{\epsilon} \\ \\
\displaystyle \qquad \qquad \qquad \qquad \qquad \qquad \qquad + \frac{E_{\epsilon,2} \, \cos\alpha - E_{\epsilon,1} \, \sin\alpha}{\sqrt{2k}} \, \D_{\alpha}f_{\epsilon} - \frac{1}{\epsilon} \, \D_{\alpha}f_{\epsilon} = 0 \, , \\ \\
\displaystyle f_{\epsilon}(\mathbf{x},k,\alpha,0) = \tilde{f}^{0}\big(x_{1}-\sqrt{2k}\,\sin\alpha, x_{2}+\sqrt{2k}\,\cos\alpha,\sqrt{2k}\,\cos\alpha, \sqrt{2k}\,\sin\alpha\big) \, , \\ \\
\displaystyle \mathbf{E}_{\epsilon}(\mathbf{x},k,\alpha,t) = \tilde{\mathbf{E}}_{\epsilon}\big(x_{1}-\sqrt{2k}\,\sin\alpha, x_{2}+\sqrt{2k}\,\cos\alpha,t\big) \, , \\ \\
\displaystyle \tilde{\mathbf{E}}_{\epsilon}(\tilde{\mathbf{x}},t) = -\nabla_{\tilde{\mathbf{x}}} \tilde{\phi}_{\epsilon}(\tilde{\mathbf{x}},t) \, , \\ \\
\displaystyle -\Delta_{\tilde{\mathbf{x}}}\tilde{\phi}_{\epsilon}(\tilde{\mathbf{x}},t) = \int_{0}^{2\pi} \hspace{-0.3cm} \int_{0}^{+\infty} \hspace{-0.3cm} f_{\epsilon}\big(\tilde{x}_{1}+\sqrt{2k}\,\sin\alpha, \tilde{x}_{2}-\sqrt{2k}\,\cos\alpha, k,\alpha, t\big) \, dk \, d\alpha - \tilde{n}_{e}(\tilde{\mathbf{x}},t) \, ,
\end{array}
\right.
\end{equation}
where $\tilde{\mathbf{x}} = (\tilde{x}_{1},\tilde{x}_{2})$ and $\mathbf{x} = (x_{1},x_{2})$,
and where $f_{\epsilon} = f_{\epsilon}(\mathbf{x},k,\alpha,t)$ and $\mathbf{E}_{\epsilon} = \mathbf{E}_{\epsilon}(\mathbf{x},k,\alpha,t)$ are linked with $\tilde{f}^{\epsilon}$ and $\tilde{\mathbf{E}}_{\epsilon}$ by 
\begin{equation} \label{def_feps_fepstilde_Eeps_Eepstilde}
\begin{split}
f_{\epsilon}(\mathbf{x},k,\alpha,t) &= \tilde{f}_{\epsilon}\big(x_{1}-\sqrt{2k}\,\sin\alpha, x_{2}+\sqrt{2k}\,\cos\alpha,\sqrt{2k}\,\cos\alpha, \sqrt{2k}\,\sin\alpha, t\big) \, , \\
\mathbf{E}_{\epsilon}(\mathbf{x},k,\alpha,t) &= \tilde{\mathbf{E}}_{\epsilon}\big(x_{1}-\sqrt{2k}\,\sin\alpha, x_{2}+\sqrt{2k}\,\cos\alpha,t\big) \, .
\end{split}
\end{equation}

\subsection{Two-scale convergence}

We set the following notations
\begin{equation}
\Omega = \R^{2} \times \R_{+} \times [0,2\pi] \, , \qquad \Gamma = \R^{2} \times \R_{+} \, , \qquad \mathcal{S} = \R_{+}\times[0,2\pi] \, ,
\end{equation}
and we consider the following Banach spaces, involving periodicity with respect to $\alpha$:
\begin{displaymath}
\begin{split}
L_{\#}^{p}\big(0,2\pi; L^{p}(\Gamma)\big) &= \Big\{ f \in L^{p}(\Omega) \, : \, 
\textit{$f$ is periodic in $\alpha$} \, \Big\} \, , \\
W_{\#}^{1,p}\big(0,2\pi; W^{1,p}(\Gamma)\big) &= \Big\{ f \in W^{1,p}(\Omega),
f(.,.,0) = f(.,.,2\pi) \, : \, 
\textit{$f$ is periodic in $\alpha$} \, \Big\} \, , \\
W_{\#}^{2,p}\big(0,2\pi; W^{2,p}(\Gamma)\big) &= \Big\{ f \in W^{2,p}(\Omega),
f(.,.,0) = f(.,.,2\pi) , \partial_\alpha f(.,.,0) = \partial_\alpha f(.,.,2\pi)  \, : 
\\ \, &~\hspace{7.4cm}
\textit{$f$ is periodic in $\alpha$} \, \Big\} \, ,
\end{split}
\end{displaymath}
and we can state  the following theorem.
\begin{theorem} \label{CV_GC}
We assume that, for a fixed $p \geq 2$, $\tilde{f}^{0}$ and $\tilde{n}_{e}$ satisfy the assumptions of Theorem \ref{CV_cart}. Then  sequences $(f_{\epsilon})_{\epsilon\,>\,0}$ and $(\tilde{\mathbf{E}}_{\epsilon})_{\epsilon\,>\,0}$
of system (\ref{VP2D_GC}) 
are bounded independently of $\epsilon$ in $L^{\infty}\big(0,T;L_{\#}^{p}(0,2\pi; L^{p}(\Gamma))\big)$ and $\big(L^{\infty}\big(0,T;W^{1,3/2}(\R^{2})\big)\big)^{2}$ respectively. As a consequence, there exist $F=F(\mathbf{x},k,\alpha,\tau,t)$ and $\tilde{\mathcal{E}} = \tilde{\mathcal{E}}(\tilde{\mathbf{x}},\tau,t)$ such that, extracting some subsequences, 
\begin{equation}
\begin{array}{rccl}
f_{\epsilon} & \longrightarrow & F & \textit{two-scale in $L^{\infty}\big(0,T; L_{\#}^{\infty}\big(0,2\pi;L_{\#}^{p}(0,2\pi; L^{p}(\Gamma))\big)\big)$,} \\
\tilde{\mathbf{E}}_{\epsilon} & \longrightarrow & \tilde{\mathcal{E}} & \textit{two-scale in $\big(L^{\infty}\big(0,T; L_{\#}^{\infty}(0,2\pi;W^{1,3/2}(\R^{2}))\big)\big)^{2}$.}
\end{array}
\end{equation}
Furthermore, there exist $G = G(\mathbf{x},k,\alpha,t) \in L^{\infty}\big(0,T;L_{\#}^{p}(0,2\pi; L^{p}(\Gamma))\big)$ and $\mathcal{E} = \mathcal{E}(\mathbf{x},k,\alpha,\tau,t) \in \big(L^{\infty}\big(0,T;L_{\#}^{\infty}(0,2\pi;W_{\#}^{1,3/2}(0,2\pi; W^{1,3/2}(\Gamma)))\big)\big)^{2}$ such that
\begin{equation} \label{relation_FG}
F(\mathbf{x},k,\alpha,\tau,t) = G(\mathbf{x},k,\alpha+\tau,t) \, ,
\end{equation}
\begin{equation} \label{def_E_Etilde}
\mathcal{E}(\mathbf{x},k,\alpha,\tau,t) = \tilde{\mathcal{E}}\big(x_{1}-\sqrt{2k}\,\sin\alpha, x_{2}+\sqrt{2k}\,\cos\alpha,\tau,t\big) \, ,
\end{equation}
and verifying
\begin{equation} \label{GC_two-scale}
\left\{
\begin{array}{l}
\displaystyle \D_{t}G + \langle \mathcal{E}_{2} \rangle \, \D_{x_{1}}G - \langle \mathcal{E}_{1} \rangle \, \D_{x_{2}}G + \langle \mathcal{F}_{k} \rangle \, \D_{k} G + \langle \mathcal{F}_{\alpha} \rangle \, \D_{\alpha} G = 0 \, , \\ \\
G(\mathbf{x},k,\alpha,0) = \cfrac{1}{2\pi} \, \tilde{f}^{0}\big(x_{1}-\sqrt{2k}\,\sin\alpha, x_{2}+\sqrt{2k}\,\cos\alpha,\sqrt{2k}\,\cos\alpha, \sqrt{2k}\,\sin\alpha\big) \, , \\ \\
\displaystyle -\nabla_{\tilde{\mathbf{x}}}\tilde{\Phi}(\tilde{\mathbf{x}},\tau,t) = \tilde{\mathcal{E}}(\tilde{\mathbf{x}},\tau,t) \, , \\ \\
\displaystyle -\Delta_{\tilde{\mathbf{x}}}\tilde{\Phi}(\tilde{\mathbf{x}},\tau,t) = \int_{\mathcal{S}} G\big(\tilde{x}_{1}+\sqrt{2k}\,\sin\alpha,\tilde{x}_{2}-\sqrt{2k}\,\cos\alpha,k,\alpha+\tau,t\big) \, dk\,d\alpha \\
\qquad \qquad \qquad \qquad \qquad \qquad \qquad \qquad \qquad \qquad \qquad \qquad \qquad \qquad - \cfrac{1}{2\pi}\,\tilde{n}_{e}(\tilde{\mathbf{x}}) \, , \\ \\
\displaystyle \mathcal{F}_{k}(\mathbf{x},k,\alpha,\tau,t) = \sqrt{2k}\,\big( \mathcal{E}_{1}(\mathbf{x},k,\alpha,\tau,t)\, \cos\alpha + \mathcal{E}_{2}(\mathbf{x},k,\alpha,\tau,t) \, \sin\alpha\big) \, , \\ \\
\displaystyle \mathcal{F}_{\alpha}(\mathbf{x},k,\alpha,\tau,t) = \frac{\mathcal{E}_{2}(\mathbf{x},k,\alpha,\tau,t) \, \cos\alpha - \mathcal{E}_{1}(\mathbf{x},k,\alpha,\tau,t) \, \sin\alpha}{\sqrt{2k}} \, ,
\end{array}
\right.
\end{equation}
where the notation $\langle \cdot \rangle$ stands for
\begin{equation} \label{bracket_def}
\langle u \rangle(\mathbf{x},k,\alpha,t) = \int_{0}^{2\pi} u(\mathbf{x},k,\alpha-\tau,\tau,t) \, d\tau \, .
\end{equation}
\end{theorem}

\begin{proof}[Proof of Theorem \ref{CV_GC}]
Several parts of this proof are only sketched since they can be redundant with \cite{Finite_Larmor_radius}. However, more details can be found in Mouton \cite{PhD_Mouton}. \\
\indent Following the same way as in \cite{Finite_Larmor_radius}, we prove that, under the assumptions of Theorem \ref{CV_cart}, we have
\begin{equation}
\big\|f_{\epsilon}(\cdot,t) \big\|_{L_{\#}^{p}(0,2\pi; L^{p}(\Gamma))} = \|\tilde{f}^{0}\|_{L^{p}(\R^{4})} \, , \forall \, t \geq 0 \, ,
\end{equation}
and, defining $\tilde{\rho}_{\epsilon}$ as
\begin{equation}
\begin{split}
\tilde{\rho}_{\epsilon}(\tilde{\mathbf{x}},t) = \int_{\mathcal{S}} f_{\epsilon}\big(\tilde{x}_{1}+\sqrt{2k}\,\sin\alpha & , \tilde{x}_{2}-\sqrt{2k}\,\cos\alpha, k, \alpha, t\big) \, dk \, d\alpha \, ,
\end{split}
\end{equation}
that the sequence $(\tilde{\rho}_{\epsilon})_{\epsilon\,>\,0}$ is bounded in $L^{\infty}\big(0,T;L^{3/2}(\R^{2})\big)$ independently of $\epsilon$. Then, we deduce that $(\tilde{\phi}_{\epsilon})_{\epsilon\,>\,0}$ and $(\tilde{\mathbf{E}}_{\epsilon})_{\epsilon\,>\,0}$ are bounded independently of $\epsilon$ in $L^{\infty}\big(0,T;W^{2,3/2}(\R^{2})\big)$ and $\big(L^{\infty}\big(0,T;W^{1,3/2}(\R^{2})\big)\big)^{2}$ respectively. As a consequence, there exist $F = F(\mathbf{x},k,\alpha,\tau,t)$, $\tilde{\Phi} = \tilde{\Phi}(\tilde{\mathbf{x}},\tau,t)$ and $\tilde{\mathcal{E}} = \tilde{\mathcal{E}}(\tilde{\mathbf{x}},\tau,t)$ such that
\begin{equation}
\begin{array}{rcll}
f_{\epsilon} & \longrightarrow & F & \textnormal{two-scale in $L^{\infty}\big(0,T; L_{\#}^{\infty}\big(0,2\pi; L_{\#}^{p}(0,2\pi; L^{p}(\Gamma))\big)\big)$} \, , \\
\tilde{\phi}_{\epsilon} & \longrightarrow & \tilde{\Phi} & \textnormal{two-scale in $L^{\infty}\big(0,T;L_{\#}^{\infty}\big(0,2\pi; W^{2,3/2}(\R^{2})\big)\big)$} \, , \\
\tilde{\mathbf{E}}_{\epsilon} & \longrightarrow & \tilde{\mathcal{E}} & \textnormal{two-scale in $\big(L^{\infty}\big(0,T;L_{\#}^{\infty}\big(0,2\pi; W^{1,3/2}(\R^{2})\big)\big)\big)^{2}$} \, .
\end{array}
\end{equation}

\indent Considering a compact set $K \subset \Gamma$, we easily remark that the sequence 
$(\phi_{\epsilon})_{\epsilon \, > \, 0}$ defined by  (\ref{def_phieps_phiepstilde}) is bounded in
$L^{\infty}\big(0,T; W_{\#}^{1,3/2}(0,2\pi; W^{1,3/2}(K))\big)$ and that all its second order derivatives
except $\D_{k}^2\phi_ \epsilon$ are bounded independently of $\epsilon$ in 
$L^{\infty}\big(0,T; L_{\#}^{3/2}\big(0,2\pi; L^{3/2}(K)))$.
The sequence
$(\mathbf{E}_{\epsilon})_{\epsilon\,>\,0}$ defined by (\ref{def_feps_fepstilde_Eeps_Eepstilde}.b) is bounded in $\big(L^{\infty}\big(0,T; W_{\#}^{1,3/2}(0,2\pi; W^{1,3/2}(K))\big)\big)^{2}$ independently of $\epsilon$. As a consequence, we claim that there exist $\Phi = \Phi(\mathbf{x},k,\alpha,\tau,t)$ and $\mathcal{E} = \mathcal{E}(\mathbf{x},k,\alpha,\tau,t)$ such that
\begin{equation} \label{CV_E_phi}
\begin{array}{rcll}
\phi_{\epsilon} & \longrightarrow & \Phi & \textnormal{two-scale in $L^{\infty}\big(0,T; L_{\#}^{\infty}\big(0,2\pi; W_{\#}^{1,3/2}(0,2\pi; W^{1,3/2}(K))\big) \big)$} \, , \\
\mathbf{E}_{\epsilon} & \longrightarrow & \mathcal{E} & \textnormal{two-scale in $\big( L^{\infty}\big(0,T; L_{\#}^{\infty}\big(0,2\pi; W_{\#}^{1,3/2}(0,2\pi; W^{1,3/2}(K))\big) \big) \big)^{2}$} \, .
\end{array}
\end{equation}
Furthermore, we remark that $\Phi$ and $\mathcal{E}$ are linked with $\tilde{\Phi}$ and $\tilde{\mathcal{E}}$ by the formula
\begin{equation}
\begin{split}
\mathcal{E}(\mathbf{x},k,\alpha,\tau,t) &= \tilde{\mathcal{E}}\big(x_{1}-\sqrt{2k}\,\sin\alpha, x_{2}+\sqrt{2k}\,\cos\alpha,\tau,t\big) \, , \\
\Phi(\mathbf{x},k,\alpha,\tau,t) &= \tilde{\Phi}\big(x_{1}-\sqrt{2k}\,\sin\alpha, x_{2}+\sqrt{2k}\,\cos\alpha,\tau,t\big) \, .
\end{split}
\end{equation}

\indent Then the vector function $\mathbf{A}_{\epsilon}$ defined by
\begin{equation}
\label{uuu123} 
\mathbf{A}_{\epsilon} = \left(
\begin{array}{c}
-\D_{x_{2}}\phi_{\epsilon} \\
\D_{x_{1}}\phi_{\epsilon} \\
\D_{\alpha}\phi_{\epsilon} \\
-\D_{k}\phi_{\epsilon}
\end{array}
\right) = \left(
\begin{array}{c}
E_{\epsilon,2} \\
-E_{\epsilon,1} \\
\sqrt{2k}\,\big( E_{\epsilon,1}\, \cos\alpha + E_{\epsilon,2} \, \sin\alpha\big) \\
\cfrac{E_{\epsilon,2} \, \cos\alpha - E_{\epsilon,1} \, \sin\alpha}{\sqrt{2k}} 
\end{array}
\right) \, , 
\end{equation}
has its three first components which are 
bounded in $L^{\infty}\big(0,T; W_{\#}^{1,3/2}(0,2\pi; W^{1,3/2}(K))\big)$,  independently of $\epsilon$
and its fourth one in $L^{\infty}\big(0,T; L_{\#}^{3/2}(0,2\pi; L^{3/2}(K))\big)$
and admits a two-scale limit denoted $\mathcal{A} = \mathcal{A}(\mathbf{x},k,\alpha,\tau,t)$ in 
$\big(L^{\infty}\big(0,T; L_{\#}^{\infty}(0,2\pi; W_{\#}^{1,3/2}(0,2\pi; W^{1,3/2}(K)))\big)\big)^{3}\times
L^{\infty}\big(0,T; L_{\#}^{\infty}(0,2\pi; L_{\#}^{3/2}(0,2\pi; L^{3/2}(K)))\big)$. 
The convergence of the three first components is the consequence of classical embedding of 
Sobolev spaces in $L^p$ spaces. Concerning the convergence of the fourth one we need to
use  that $E_{\epsilon,1}$ and $E_{\epsilon,2}$ are bounded independently of $\epsilon$ in  
$L^{\infty}\big(0,T; W_{\#}^{1,3/2}(0,2\pi;$ $ W^{1,3/2}(K))\big)$ and consequently in 
$L^{\infty}\big(0,T; L_{\#}^{6}(0,2\pi; L^{6}(K))\big)$ and that $(2k)^{-1/2}$ is, for any $k_{\max}$ 
and any $q<2$, in $L^{q}(0,k_{\max})$. The vector function $\mathcal{A}$ is linked with $\Phi$ and $\mathcal{E}$ as follows:
\begin{equation}
\mathcal{A} = \left(
\begin{array}{c}
-\D_{x_{2}}\Phi \\
\D_{x_{1}}\Phi \\
\D_{\alpha}\Phi \\
-\D_{k}\Phi
\end{array}
\right) = \left(
\begin{array}{c}
\mathcal{E}_{2} \\
-\mathcal{E}_{1} \\
\sqrt{2k}\,\big( \mathcal{E}_{1}\, \cos\alpha + \mathcal{E}_{2} \, \sin\alpha\big) \\
\cfrac{\mathcal{E}_{2} \, \cos\alpha - \mathcal{E}_{1} \, \sin\alpha}{\sqrt{2k}} 
\end{array}
\right) \, .
\end{equation}

In order to establish the two-scale limit model, we cannot simply apply Theorem 1.3 of \cite{Finite_Larmor_radius}: indeed, the formulation (\ref{VP2D_GC}.a) of Vlasov equation does not fit with the assumptions which are needed for applying this theorem since the differential operator $f \mapsto -\cfrac{1}{\epsilon}\,\D_{\alpha}f$ cannot be written under the form
\begin{equation}
f \mapsto \cfrac{1}{\epsilon} \, \Big(\mathbb{M}\,\left(
\begin{array}{c}
x_{1} \\ x_{2} \\ k \\ \alpha
\end{array}
\right) + \mathbf{N} \Big) \cdot \nabla f \, ,
\end{equation}
where $\mathbb{M}$ is a constant square matrix satisfying $Tr(\mathbb{M}) = 0$, and $\mathbf{N} \in Im(\mathbb{M})$. However, the approach which is considered in \cite{Finite_Larmor_radius} can be adapted to the present case. \\

\indent Firstly, we prove that there exists a function $G$ such that $F(\mathbf{x},k,\alpha,\tau,t) = G(\mathbf{x},k,\alpha+\tau,t)$. To reach such a result, we consider a test function $\psi = \psi(\mathbf{x},k,\alpha,\tau,t)$ on $\Omega \times [0,2\pi] \times [0,T]$ which is $2\pi$-periodic in $\alpha$ and $\tau$. If we multiply (\ref{VP2D_GC}.a) by $\psi(\mathbf{x},k,\alpha,\frac{t}{\epsilon},t)$ and integrate over $\Omega \times [0,T]$, we obtain
\begin{equation} \label{weak_VGC}
\begin{split}
\int_{0}^{T} \int_{\Omega} f_{\epsilon}&(\mathbf{x},k,\alpha,t) \, \Big[ \D_{t}\psi\big(\mathbf{x},k,\alpha,\frac{t}{\epsilon},t\big) + \frac{1}{\epsilon}\,\D_{\tau}\psi\big(\mathbf{x},k,\alpha,\frac{t}{\epsilon},t\big) \\
&\quad + \mathbf{A}_{\epsilon}(\mathbf{x},k,\alpha,t) \cdot \nabla f_{\epsilon}(\mathbf{x},k,\alpha,t) - \frac{1}{\epsilon}\,\D_{\alpha}\psi\big(\mathbf{x},k,\alpha,\frac{t}{\epsilon},t\big) \Big] \, d\mathbf{x} \, dk \, d\alpha \, dt \\
& = -\int_{\Omega} \tilde{f}^{0}(x_{1}-\sqrt{2k}\,\sin\alpha,x_{2}+\sqrt{2k}\,\cos\alpha,\sqrt{2k}\,\cos\alpha,\sqrt{2k}\,\sin\alpha) \\
&\qquad \qquad \qquad \times \psi(x_{1},x_{2},k,\alpha,0,0) \, dx_{1}\, dx_{2} \, dk \, d\alpha \, .
\end{split}
\end{equation}
Multiplying (\ref{weak_VGC}) by $\epsilon$ and letting $\epsilon \to 0$, we obtain the weak formulation of $\D_{\tau}F - \D_{\alpha}F = 0$, which indicates that there exists a function $G \in L^{\infty}\big(0,T;L_{\#}^{p}(0,2\pi; L^{p}(\Gamma))\big)$ such that
\begin{equation}
F(\mathbf{x},k,\alpha,\tau,t) = G(\mathbf{x},k,\alpha+\tau,t) \, .
\end{equation}
\indent Secondly, we introduce the sequence $(g_{\epsilon})_{\epsilon\,>\,0}$ defined by
\begin{equation}
g_{\epsilon}(\mathbf{x},k,\alpha,t) = f_{\epsilon}\Big(\mathbf{x},k,\alpha-\cfrac{t}{\epsilon},t\Big) \, .
\end{equation}
In the spirit of \cite{Finite_Larmor_radius}, we prove that $g_{\epsilon}$ strongly converges to $2\pi\,G$ in a given Banach space. For that, we notice that, up to a subsequence, $g_{\epsilon}$ two-scale converges to $G$ in $L^{\infty}\big(0,T;L_{\#}^{\infty}\big(0,2\pi;$ $ L_{\#}^{p}(0,2\pi; L^{p}(\Gamma))\big)\big)$ since we have
\begin{equation}
\begin{split}
&\int_{0}^{T} \int_{\Omega} g_{\epsilon}(\mathbf{x},k,\alpha,t) \, \psi\big(\mathbf{x},k,\alpha,\frac{t}{\epsilon},t\big) \, d\mathbf{x} \ dk \, d\alpha \, dt \\
&\qquad = \int_{0}^{T} \int_{\Omega} f_{\epsilon}(\mathbf{x},k,\alpha,t) \, \psi\big(\mathbf{x},k,\alpha+\frac{t}{\epsilon},\frac{t}{\epsilon},t\big) \, d\mathbf{x} \ dk \, d\alpha \, dt \\
&\qquad \to \int_{0}^{2\pi} \int_{0}^{T} \int_{\Omega} F(\mathbf{x},k,\alpha,\tau,t) \, \psi(\mathbf{x},k,\alpha+\tau,\tau,t) \, d\mathbf{x} \ dk \, d\alpha \, dt \, d\tau \\
&\qquad = \int_{0}^{2\pi} \int_{0}^{T} \int_{\Omega} G(\mathbf{x},k,\alpha,t) \, \psi(\mathbf{x},k,\alpha,\tau,t) \, d\mathbf{x} \ dk \, d\alpha \, dt \, d\tau \, ,
\end{split}
\end{equation}
for any test function $\psi$ on $\Omega \times [0,2\pi] \times [0,T]$ which is $2\pi$-periodic in $\alpha$ and $\tau$. As a consequence, $g_{\epsilon}$ weakly-* converges to $2\pi\,G$ in $L^{\infty}\big(0,T;L_{\#}^{p}(0,2\pi; L^{p}(\Gamma))\big)$. \\
\indent Let us prove that this weak-* convergence is a strong convergence in a given Banach space. This is the aim of the following lemma:

\begin{lemma} \label{CV_strong_compact}
For any compact subset $K$ of $\Gamma$, and up to a subsequence, $g_{\epsilon}$ strongly converges to $2\pi G$ in $L^{\infty}\big(0,T;(W_{\#}^{1,3/2}(0,2\pi; W_{0}^{1,3/2}(K)))^{*}\big)$.
\end{lemma}
\noindent
In this Lemma $(W_{\#}^{1,3/2}(0,2\pi; W_{0}^{1,3/2}(K)))^{*}$ stands for the dual of
$W_{\#}^{1,3/2}(0,2\pi; W_{0}^{1,3/2}(K))$.

\begin{proof}[Proof of Lemma \ref{CV_strong_compact}]
For any compact subset $K$ of $\Gamma$, $(g_{\epsilon})_{\epsilon\,>\,0}$ and $(\mathbf{A}_{\epsilon})_{\epsilon\,>\,0}$ are bounded independently of $\epsilon$ in $L^{\infty}\big(0,T;L_{\#}^{p}(0,2\pi; L^{p}(K))\big)$ and 
$\big(L^{\infty}\big(0,T; W_{\#}^{1,3/2}(0,2\pi; W^{1,3/2}(K))\big)\big)^{3}$ 
$\times L^{\infty}\big(0,T; L_{\#}^{3/2}(0,2\pi; L^{3/2}(K))\big)$ respectively.
 Then, remarking that $g_{\epsilon}$ is solution of
\begin{equation} \label{eq_g_GC}
\left\{
\begin{array}{l}
\displaystyle \D_{t}g_{\epsilon}(\mathbf{x},k,\alpha,t) + \mathbf{A}_{\epsilon}\big(\mathbf{x},k,\alpha-\frac{t}{\epsilon},t\big) \cdot \nabla g_{\epsilon}(\mathbf{x},k,\alpha,t) = 0 \, , \\ \\
g_{\epsilon}(\mathbf{x},k,\alpha,0) = \tilde{f}^{0}(x_{1}-\sqrt{2k}\,\sin\alpha,x_{2}+\sqrt{2k}\,\cos\alpha,\sqrt{2k}\,\cos\alpha,\sqrt{2k}\,\sin\alpha) \, ,
\end{array}
\right.
\end{equation}
we use similar arguments as the ones given after equation (\ref {uuu123}) to deduce
\begin{enumerate}
\item $(\mathbf{A}_{\epsilon})_{\epsilon\,>\,0}$ is bounded independently of $\epsilon$ in $\big(L^{\infty}\big(0,T;L_{\#}^{3/2}(0,2\pi; L^{3/2}(K))\big)\big)^{4}$,
\item $\big(g_{\epsilon}(\mathbf{x},k,\alpha,t)\,\mathbf{A}_{\epsilon}(\mathbf{x},k,\alpha-\frac{t}{\epsilon},t)\big)_{\epsilon\,>\,0}$ is bounded in $\big(L^{\infty}\big(0,T;L_{\#}^{r}(0,2\pi; L^{r}(K))\big)\big)^{4}$ independently of $\epsilon$ with $r$ defined by $\frac{1}{r} = \frac{1}{p}+\frac{1}{q}$ 
($r \in \, ]1,\frac{3}{2}[$),
\item $(\D_{t}g_{\epsilon})_{\epsilon\,>\,0}$ is bounded in $L^{\infty}\big(0,T;\big(W_{\#}^{1,r^{*}}(0,2\pi; W_{0}^{1,r^{*}}(K))\big)^{*}\big)$ independently of $\epsilon$ with $\frac{1}{r^{*}}+\frac{1}{r} = 1$.
\end{enumerate}
If {$r < 3/2$}, the embedding $W_{\#}^{1,r^{*}}(0,2\pi; W_{0}^{1,r^{*}}(K)) \subset W_{\#}^{1,3/2}(0,2\pi; W_{0}^{1,3/2}(K))$ is compact with density. Furthermore, Rellich-Kondrakov's theorem (see \cite{Adams}) gives the compact embedding $L_{\#}^{p}(0,2\pi; L^{p}(K)) \subset (W_{\#}^{1,3/2}(0,2\pi; W_{0}^{1,3/2}(K)))^{*}$ since $p \geq 2$. Then, we apply Aubin-Lions' lemma (see \cite{Lions}) and we prove that the functional space $\mathcal{U}$ defined by
\begin{equation}
\begin{array}{c}
\begin{split}
\displaystyle \mathcal{U} = \Big\{ u \in L^{\infty}\big(0,T;L_{\#}^{p}(0,2\pi; L^{p}&(K))\big) \, : \\
& \D_{t}u \in L^{\infty}\big(0,T;\big(W_{\#}^{1,r^{*}}(0,2\pi; W_{0}^{1,r^{*}}(K))\big)^{*}\big) \Big\} \, , 
\end{split}
\end{array}
\end{equation}
is compactly embedded in $L^{\infty}\big(0,T;(W_{\#}^{1,3/2}(0,2\pi; W_{0}^{1,3/2}(K)))^{*}\big)$. Since $g_{\epsilon} \in \mathcal{U}$ for all $\epsilon$, we deduce that the weak-* convergence of $g_{\epsilon}$ to $2\pi\,G$ in $L^{\infty}\big(0,T;L_{\#}^{p}(0,2\pi; L^{p}(K))\big)$ is a strong convergence in $L^{\infty}\big(0,T; (W_{\#}^{1,3/2}(0,2\pi; W_{0}^{3/2}(K)))^{*}\big)$. \\
If {$r \geq 3/2$}, the compact embedding $\big(W_{\#}^{1,r^{*}}(0,2\pi; W_{0}^{1,r^{*}}(K))\big)^{*}\subset \big(W_{\#}^{1,3/2}(0,2\pi; W_{0}^{1,3/2}(K))\big)^{*}$ is gotten directly. If we introduce the functional space $\mathcal{U}'$ defined by
\begin{equation}
\begin{array}{c}
\begin{split}
\displaystyle \mathcal{U}' = \Big\{ u \in L^{\infty}\big(0,T;L_{\#}^{p}(0,2\pi; L^{p}&(K))\big) \, : \\
& \D_{t}u \in L^{\infty}\big(0,T;\big(W_{\#}^{1,3/2}(0,2\pi; W_{0}^{1,3/2}(K))\big)^{*}\big) \Big\} \, , 
\end{split}
\end{array}
\end{equation}
we remark that the sequence $(g_{\epsilon})_{\epsilon\,>\,0}$ is bounded in $\mathcal{U}'$ independently of $\epsilon$. By using Aubin-Lions' lemma, we prove that $\mathcal{U}' \subset L^{\infty}\big(0,T;(W_{\#}^{1,3/2}(0,2\pi; W_{0}^{1,3/2}(K)))^{*}\big)$ is a compact embedding. Then, $g_{\epsilon}$ strongly converges to $2\pi \,G$ in $L^{\infty}\big(0,T; (W_{\#}^{1,3/2}(0,2\pi; W_{0}^{1,3/2}(K)))^{*}\big)$.
\end{proof}

\textit{ \\ }

Let us finish the proof of Theorem \ref{CV_GC} by establishing a transport equation satisfied by $G$. Let us consider a test function $\psi = \psi(\mathbf{x},k,\alpha,t)$ on $\Omega$. We denote its compact support in $\Gamma$ by $K$ and we assume that $\psi$ and its first order derivatives are $2\pi$-periodic in the $\alpha$ direction. Then we have
\begin{equation} \label{weak_form_geps}
\begin{split}
\int_{0}^{T} \int_{0}^{2\pi} & \int_{K} g_{\epsilon}(\mathbf{x},k,\alpha,t) \, \D_{t}\psi(\mathbf{x},k,\alpha,t) d\mathbf{x} \, dk \, d\alpha \, dt \\
& + \int_{0}^{T} \int_{0}^{2\pi} \, \int_{K} g_{\epsilon}(\mathbf{x},k,\alpha,t) \, \mathbf{A}_{\epsilon}\big(\mathbf{x},k,\alpha-\frac{t}{\epsilon},t\big) \cdot \nabla \psi(\mathbf{x},k,\alpha,t) \Big] \, d\mathbf{x} \, dk \, d\alpha \, dt \\
& + \int_{0}^{2\pi} \int_{K} \tilde{f}^{0}(x_{1}-\sqrt{2k}\,\sin\alpha,x_{2}+\sqrt{2k}\,\cos\alpha,\sqrt{2k}\,\cos\alpha,\sqrt{2k}\,\sin\alpha) \\
&\qquad \qquad \qquad \qquad \qquad \qquad \qquad \qquad \qquad \qquad \times \psi(\mathbf{x},k,\alpha,0) \, d\mathbf{x} \, dk \, d\alpha = 0 \, .
\end{split}
\end{equation}
Since the vector function $\mathbf{Z}_{\epsilon}(\mathbf{x},k,\alpha,t) = g_{\epsilon}(\mathbf{x},k,\alpha,t) \, \mathbf{A}_{\epsilon}(\mathbf{x},k,\alpha-\frac{t}{\epsilon},t)$ is bounded independently of $\epsilon$ in $\big(L^{\infty}\big(0,T;L_{\#}^{r}(0,2\pi; L^{r}(K))\big)\big)^{4}$ with $r \in \, ]1,\frac{3}{2}[$, there exists a function $\mathcal{Z} = \mathcal{Z}(\mathbf{x},k,\alpha,\tau,t)$ in $\big(L^{\infty}\big(0,T;L_{\#}^{\infty}(0,2\pi;L_{\#}^{r}(0,2\pi; L^{r}(K)))\big)\big)^{4}$ such that $\mathbf{Z}_{\epsilon}$ two-scale converges to $\mathcal{Z}$ in $\big(L^{\infty}\big(0,T;L_{\#}^{\infty}(0,2\pi;L_{\#}^{r}(0,2\pi; L^{r}(K)))\big)\big)^{4}$ and weakly-* converges to $\displaystyle \int_{0}^{2\pi}\mathcal{Z}(\cdot,\cdot,\cdot,\tau,\cdot)\,d\tau$ in $\big(L^{\infty}\big(0,T;L_{\#}^{r}(0,2\pi; L^{r}(K))\big)\big)^{4}$. \\
\indent If we now restrict our choice of test functions $\psi$ to the one having their support satisfying
\begin{equation}
K \, \cap \, \big(\R^{2} \times \{0\}\big) = \emptyset \, ,
\end{equation}
then, $\mathbf{A}_{\epsilon}(\mathbf{x},k,\alpha-\frac{t}{\epsilon},t)$ is bounded in $\big(L^{\infty}\big(0,T;W_{\#}^{1,3/2}(0,2\pi; W^{1,3/2}(K))\big)\big)^{4}$ independently of $\epsilon$ and admits a two-scale limit in space $\big(L^{\infty}\big(0,T;$ $L_{\#}^{\infty}(0,2\pi;W_{\#}^{1,3/2}(0,2\pi; $ $W^{1,3/2}(K)))\big)\big)^{4}$ which is exactly $\mathcal{A}$. Since sequence $(g_{\epsilon})_{\epsilon\,>\,0}$ converges to $2\pi\,G$ in space $L^{\infty}\big(0,T;(W_{\#}^{1,3/2}(0,2\pi; $ $W_{0}^{1,3/2}(K)))^{*}\big)$ strongly, we can deduce that
\begin{equation}
\mathcal{Z} = 2\pi\,G\,\mathcal{A} \, ,
\end{equation}
almost everywhere on $\Omega\times[0,T]\setminus\{(\mathbf{x},k,\alpha,t),k=0\}$ and then
almost everywhere on $\Omega\times[0,T]$, and we obtain the weak formulation of
\begin{equation}
\left\{
\begin{array}{l}
\displaystyle \D_{t}G(\mathbf{x},k,\alpha,t) + \Bigg[\int_{0}^{2\pi} \mathcal{A}(\mathbf{x},k,\alpha-\tau,\tau,t) \,d\tau \Bigg] \cdot \nabla G(\mathbf{x},k,\alpha,t) = 0 \, , \\
\displaystyle G(\mathbf{x},k,\alpha,0) = \cfrac{1}{2\pi} \, \tilde{f}^{0}(x_{1}-\sqrt{2k}\,\sin\alpha,x_{2}+\sqrt{2k}\,\cos\alpha,\sqrt{2k}\,\cos\alpha,\sqrt{2k}\,\sin\alpha) \, ,
\end{array}
\right.
\end{equation}
when $\epsilon \to 0$ in (\ref{weak_form_geps}). Finally, we extend this conclusion to every compact subsets $K$ of $\Gamma$ by remarking that the measure of $\R^{2} \times \{0\}$ is null. \\

\indent In order to obtain Poisson type equations (\ref{GC_two-scale}.c) and (\ref{GC_two-scale}.d), we consider a test function $\tilde{\psi} = \tilde{\psi}(\tilde{\mathbf{x}},\tau,t)$ on $\R^{2} \times [0,2\pi] \times [0,T]$ which is $2\pi$-periodic in $\tau$, we multiply $\tilde{\mathbf{E}}_{\epsilon}(\tilde{\mathbf{x}},t)$, $\nabla_{\tilde{\mathbf{x}}}\tilde{\phi}_{\epsilon}(\tilde{\mathbf{x}},t)$ and $\Delta_{\tilde{\mathbf{x}}}\tilde{\phi}_{\epsilon}(\tilde{\mathbf{x}},t)$ by $\tilde{\psi}(\tilde{\mathbf{x}},\frac{t}{\epsilon},t)$ and we integrate in $\tilde{\mathbf{x}}$ and $t$. We obtain
\begin{equation} \label{two-scale_Poisson_1}
\int_{0}^{T} \int_{\R^{2}} \tilde{\mathbf{E}}_{\epsilon}(\tilde{\mathbf{x}},t) \, \psi\big(\tilde{\mathbf{x}},\frac{t}{\epsilon},t\big) \, d\tilde{\mathbf{x}}\,dt \to \int_{0}^{2\pi} \int_{0}^{T} \int_{\R^{2}} \tilde{\mathcal{E}}(\tilde{\mathbf{x}},\tau,t) \, \psi(\tilde{\mathbf{x}},\tau,t) \, d\tilde{\mathbf{x}}\,dt \, d\tau \, ,
\end{equation}
\begin{equation} \label{two-scale_Poisson_2}
\int_{0}^{T} \int_{\R^{2}} \nabla_{\tilde{\mathbf{x}}}\tilde{\phi}_{\epsilon}(\tilde{\mathbf{x}},t) \, \psi\big(\tilde{\mathbf{x}},\frac{t}{\epsilon},t\big) \, d\tilde{\mathbf{x}}\,dt \to \int_{0}^{2\pi} \int_{0}^{T} \int_{\R^{2}} \nabla_{\tilde{\mathbf{x}}} \tilde{\Phi}(\tilde{\mathbf{x}},\tau,t) \, \psi(\tilde{\mathbf{x}},\tau,t) \, d\tilde{\mathbf{x}}\,dt \, d\tau \, ,
\end{equation}
\begin{equation} \label{two-scale_Poisson_3}
\int_{0}^{T} \int_{\R^{2}} \Delta_{\tilde{\mathbf{x}}}\tilde{\phi}_{\epsilon}(\tilde{\mathbf{x}},t) \, \psi\big(\tilde{\mathbf{x}},\frac{t}{\epsilon},t\big) \, d\tilde{\mathbf{x}}\,dt \to \int_{0}^{2\pi} \int_{0}^{T} \int_{\R^{2}} \Delta_{\tilde{\mathbf{x}}} \tilde{\Phi}(\tilde{\mathbf{x}},\tau,t) \, \psi(\tilde{\mathbf{x}},\tau,t) \, d\tilde{\mathbf{x}}\,dt \, d\tau \, ,
\end{equation}
when $\epsilon$ converges to 0. Since $F$ is the two-scale limit of $(f_{\epsilon})_{\epsilon\,>\,0}$, we also have
\begin{equation} \label{two-scale_Poisson_4}
\begin{split}
&\int_{0}^{T} \int_{\R^{2}} \Bigg(\int_{\mathcal{S}} f_{\epsilon}\big(\tilde{x}_{1}+\sqrt{2k}\,\sin\alpha,\tilde{x}_{2}-\sqrt{2k}\,\cos\alpha,k,\alpha,t) \, dk\,d\alpha \Bigg) \, \psi\big(\tilde{\mathbf{x}},\frac{t}{\epsilon},t\big) \, d\tilde{\mathbf{x}}\,dt \\
& \to \int_{0}^{2\pi} \int_{0}^{T} \int_{\R^{2}} \int_{\mathcal{S}} F\big(\tilde{x}_{1}+\sqrt{2k}\,\sin\alpha,\tilde{x}_{2}-\sqrt{2k}\,\cos\alpha,k,\alpha,\tau,t) \\ 
& \qquad \qquad \qquad \qquad \qquad \qquad \qquad \qquad \qquad \qquad \qquad \times \psi(\tilde{\mathbf{x}},\tau,t) \, dk\,d\alpha \, d\tilde{\mathbf{x}}\,dt \, d\tau \, ,
\end{split}
\end{equation}
when $\epsilon \to 0$. Then, gathering convergence results (\ref{two-scale_Poisson_1})-(\ref{two-scale_Poisson_4}), we obtain the weak formulation of (\ref{GC_two-scale}.c) and (\ref{GC_two-scale}.d). 
\end{proof}

\subsection{Weak-* convergence under non-physical hypothesis}

\indent Firstly, we have a direct corollary of Theorem \ref{CV_GC}:
\begin{corollary}
Up to some subsequences, 
\begin{itemize}
\item $f_{\epsilon}$ weakly-* converges to $f \in L^{\infty}\big(0,T;L_{\#}^{p}(0,2\pi; L^{p}(\Gamma))\big)$ with $f$ defined by
\begin{equation} \label{link_f_F}
f(\mathbf{x},k,\alpha,t) = \int_{0}^{2\pi}F(\mathbf{x},k,\alpha,\tau,t) \, d\tau \, ,
\end{equation}
\item $\tilde{\mathbf{E}}_{\epsilon}$ weakly-* converges to $\tilde{\mathbf{E}} \in \big(L^{\infty}\big(0,T;W^{1,3/2}(\R^{2})\big)\big)^{2}$ with $\tilde{\mathbf{E}}$ defined by
\begin{equation}
\tilde{\mathbf{E}}(\mathbf{x},t) = \int_{0}^{2\pi}\tilde{\mathcal{E}}(\mathbf{x},\tau,t) \, d\tau \, ,
\end{equation}
\item $\tilde{\phi}_{\epsilon}$ weakly-* converges to $\tilde{\phi} \in L^{\infty}\big(0,T;W^{2,3/2}(\R^{2})\big)$ with $\tilde{\phi}$ defined by
\begin{equation}
\tilde{\phi}(\mathbf{x},t) = \int_{0}^{2\pi}\tilde{\Phi}(\mathbf{x},\tau,t) \, d\tau \, .
\end{equation}
\end{itemize}
\end{corollary}

From now, we add the following non-physical assumption:
\begin{equation} \label{hypothesis_E}
\tilde{\mathbf{E}}_{\epsilon} \to \tilde{\mathbf{E}} \, \, \, \, \textnormal{strongly.}
\end{equation}

\begin{lemma} \label{Phi_property}
Under the hypotheses of Theorem \ref{CV_GC} and (\ref{hypothesis_E}), we have
\begin{equation}
\forall\, (\mathbf{x},k,\alpha,t) \in \Omega \times [0,T] \, , \qquad \D_{\alpha} \Bigg( \int_{0}^{2\pi} \Phi(\mathbf{x},k,\alpha-\tau,\tau,t) \, d\tau \Bigg) = 0 \, .
\end{equation}
\end{lemma}

\begin{proof}[Proof of Lemma \ref{Phi_property}]
We consider a compact subset $K$ of $\Gamma$ and a test function $\psi = \psi(\mathbf{x},k,\alpha,t)$ on $\Omega \times [0,T]$ which support in $(\mathbf{x},k)$ is included in $K$. Since $(\phi_{\epsilon})_{\epsilon\,>\,0}$ is bounded independently of $\epsilon$ is $L^{\infty}\big(0,T;W_{\#}^{1,3/2}(0,2\pi; W^{1,3/2}(K))\big)$ and that all its second order derivatives, except $\D_{k}^{2}\phi_{\epsilon}$, are also bounded in $L^{\infty}\big(0,T;L_{\#}^{3/2}(0,2\pi; L^{3/2}(K))\big)$, we have
\begin{equation}
\begin{split}
\int_{0}^{T}&\int_{0}^{2\pi}\int_{\Gamma} \nabla_{\mathbf{x}} \D_{\alpha}\phi_{\epsilon}(\mathbf{x},k,\alpha-\cfrac{t}{\epsilon},t) \, \psi(\mathbf{x},k,\alpha,t) \, d\mathbf{x} \, dk \, d\alpha \, dt \\
&= - \int_{0}^{T}\int_{0}^{2\pi}\int_{\Gamma} \nabla_{\mathbf{x}} \phi_{\epsilon}(\mathbf{x},k,\alpha,t) \, (\D_{\alpha}\psi)(\mathbf{x},k,\alpha+\cfrac{t}{\epsilon},t) \, d\mathbf{x} \, dk \, d\alpha \, dt \\
&= \int_{0}^{T}\int_{0}^{2\pi}\int_{\Gamma} \mathbf{E}_{\epsilon}(\mathbf{x},k,\alpha,t) \, (\D_{\alpha}\psi)(\mathbf{x},k,\alpha+\cfrac{t}{\epsilon},t) \, d\mathbf{x} \, dk \, d\alpha \, dt \\
&= \int_{0}^{T}\int_{0}^{2\pi}\int_{\Gamma} \tilde{\mathbf{E}}_{\epsilon}(\tilde{\mathbf{x}},t) \\
&\qquad \qquad \times (\D_{\alpha}\psi)\big(\tilde{x}_{1}+\sqrt{2k}\,\sin\alpha,\tilde{x}_{2}-\sqrt{2k}\,\cos\alpha,k,\alpha+\cfrac{t}{\epsilon},t\big) \, d\tilde{\mathbf{x}} \, dk \, d\alpha \, dt \, .
\end{split}
\end{equation}
Since $\tilde{\mathbf{E}}_{\epsilon}$ converges to $\tilde{\mathbf{E}}$ strongly, we obtain
\begin{equation}
\begin{split}
&\lim_{\epsilon\,\to\,0} \int_{0}^{T}\int_{0}^{2\pi}\int_{\Gamma} \nabla_{\mathbf{x}} \D_{\alpha}\phi_{\epsilon}(\mathbf{x},k,\alpha-\cfrac{t}{\epsilon},t) \, \psi(\mathbf{x},k,\alpha,t) \, d\mathbf{x} \, dk \, d\alpha \, dt \\
& \quad = \int_{0}^{T}\int_{0}^{2\pi}\int_{\Gamma} \tilde{\mathbf{E}}(\tilde{\mathbf{x}},t) \\
&\qquad \qquad \times \int_{0}^{2\pi} (\D_{\alpha}\psi)\big(\tilde{x}_{1}+\sqrt{2k}\,\sin\alpha,\tilde{x}_{2}-\sqrt{2k}\,\cos\alpha,k,\alpha+\tau,t\big) \, d\tau \, d\tilde{\mathbf{x}} \, dk \, d\alpha \, dt \\
& \quad = \int_{0}^{T}\int_{0}^{2\pi}\int_{\Gamma} \tilde{\mathbf{E}}(\tilde{\mathbf{x}},t) \\
&\qquad \qquad \times \int_{0}^{2\pi} \D_{\tau}\Big(\psi\big(\tilde{x}_{1}+\sqrt{2k}\,\sin\alpha,\tilde{x}_{2}-\sqrt{2k}\,\cos\alpha,k,\alpha+\tau,t\big)\Big) \, d\tau \, d\tilde{\mathbf{x}} \, dk \, d\alpha \, dt \\
&\quad = 0 \, .
\end{split}
\end{equation}
On another hand, we deduce from (\ref{CV_E_phi}) that we also have
\begin{equation}
\begin{split}
\int_{0}^{T}&\int_{0}^{2\pi}\int_{\Gamma} \nabla_{\mathbf{x}} \D_{\alpha}\phi_{\epsilon}(\mathbf{x},k,\alpha-\cfrac{t}{\epsilon},t) \, \psi(\mathbf{x},k,\alpha,t) \, d\mathbf{x} \, dk \, d\alpha \, dt \\
& \to \int_{0}^{T}\int_{0}^{2\pi}\int_{\Gamma} \nabla_{\mathbf{x}} \D_{\alpha} \Bigg(\int_{0}^{2\pi} \Phi(\mathbf{x},k,\alpha-\tau,\tau,t) \,d\tau \Bigg)\, \psi(\mathbf{x},k,\alpha,t) \, d\mathbf{x} \, dk \, d\alpha \, dt \, .
\end{split}
\end{equation}
Then it is straightforward that 
\begin{equation}
\nabla_{\mathbf{x}} \D_{\alpha} \Bigg(\int_{0}^{2\pi} \Phi(\mathbf{x},k,\alpha-\tau,\tau,t) \,d\tau \Bigg) = 0 \, ,
\end{equation}
in $L^{\infty}\big(0,T;L_{\#}^{3/2}(0,2\pi; L^{3/2}(K))\big)$ and that $\displaystyle (\mathbf{x},k,\alpha,t) \mapsto \D_{\alpha} \Bigg(\int_{0}^{2\pi} \Phi(\mathbf{x},k,\alpha-\tau,\tau,t) \,d\tau \Bigg)$ does not depend on $\mathbf{x}$. \\
\indent Following the same approach, we use the fact that the sequences $(\D_{k,\alpha}^{2}\phi_{\epsilon})_{\epsilon\,>\,0}$ and $(\D_{\alpha}^{2}\phi_{\epsilon})_{\epsilon\,>\,0}$ are bounded in $L^{\infty}\big(0,T;L_{\#}^{3/2}(0,2\pi;L^{3/2}(K))\big)$ independently of $\epsilon$ to obtain
\begin{equation}
\begin{split}
\D_{k,\alpha}^{2} \phi_{\epsilon}(\mathbf{x},k,\alpha-\cfrac{t}{\epsilon},t) & \qquad \stackrel{*}{\rightharpoonup} \qquad \D_{k,\alpha}^{2} \Bigg(\int_{0}^{2\pi} \Phi(\mathbf{x},k,\alpha-\tau,\tau,t) \,d\tau \Bigg) \, , \\
\D_{\alpha}^{2} \phi_{\epsilon}(\mathbf{x},k,\alpha-\cfrac{t}{\epsilon},t) & \qquad \stackrel{*}{\rightharpoonup} \qquad \D_{\alpha}^{2}\Bigg(\int_{0}^{2\pi} \Phi(\mathbf{x},k,\alpha-\tau,\tau,t) \,d\tau \Bigg) \, , 
\end{split}
\end{equation}
in $L^{\infty}\big(0,T;L_{\#}^{3/2}(0,2\pi;L^{3/2}(K))\big)$. On the other hand, we have
\begin{equation}
\begin{split}
&\int_{0}^{T}\int_{0}^{2\pi}\int_{\Gamma} \D_{k,\alpha}^{2}\phi_{\epsilon}(\mathbf{x},k,\alpha-\cfrac{t}{\epsilon},t) \, \psi(\mathbf{x},k,\alpha,t) \, d\mathbf{x} \, dk \, d\alpha \, dt \\
&\quad = \int_{0}^{T}\int_{0}^{2\pi}\int_{\Gamma}  \cfrac{E_{\epsilon,2}(\mathbf{x},k,\alpha,t)\,\cos\alpha}{\sqrt{2k}}\, (\D_{\alpha}\psi)(\mathbf{x},k,\alpha+\cfrac{t}{\epsilon},t) \, d\mathbf{x} \, dk \, d\alpha \, dt \\
&\qquad - \int_{0}^{T}\int_{0}^{2\pi}\int_{\Gamma}  \cfrac{E_{\epsilon,1}(\mathbf{x},k,\alpha,t)\,\sin\alpha}{\sqrt{2k}}\, (\D_{\alpha}\psi)(\mathbf{x},k,\alpha+\cfrac{t}{\epsilon},t) \, d\mathbf{x} \, dk \, d\alpha \, dt \\
&\quad = \int_{0}^{T}\int_{0}^{2\pi}\int_{\Gamma} \tilde{E}_{\epsilon,2}(\tilde{\mathbf{x}},t)\, \cfrac{\cos\alpha}{\sqrt{2k}} \\
&\quad \qquad \qquad \times (\D_{\alpha}\psi)(\tilde{x}_{1}+\sqrt{2k}\,\sin\alpha, \tilde{x}_{2}-\sqrt{2k}\,\cos\alpha,k,\alpha+\cfrac{t}{\epsilon},t) \, d\tilde{\mathbf{x}} \, dk \, d\alpha \, dt  \\
&\quad \qquad - \int_{0}^{T}\int_{0}^{2\pi}\int_{\Gamma} \tilde{E}_{\epsilon,1}(\tilde{\mathbf{x}},t)\, \cfrac{\sin\alpha}{\sqrt{2k}} \\
&\qquad \qquad \times (\D_{\alpha}\psi)(\tilde{x}_{1}+\sqrt{2k}\,\sin\alpha, \tilde{x}_{2}-\sqrt{2k}\,\cos\alpha,k,\alpha+\cfrac{t}{\epsilon},t) \, d\tilde{\mathbf{x}} \, dk \, d\alpha \, dt \, ,
\end{split}
\end{equation}
and
\begin{equation}
\begin{split}
&\int_{0}^{T}\int_{0}^{2\pi}\int_{\Gamma} \D_{\alpha}^{2}\phi_{\epsilon}(\mathbf{x},k,\alpha-\cfrac{t}{\epsilon},t) \, \psi(\mathbf{x},k,\alpha,t) \, d\mathbf{x} \, dk \, d\alpha \, dt \\
&\quad = \int_{0}^{T}\int_{0}^{2\pi}\int_{\Gamma} \sqrt{2k}\,E_{\epsilon,1}(\mathbf{x},k,\alpha,t)\,\cos\alpha\, (\D_{\alpha}\psi)(\mathbf{x},k,\alpha+\cfrac{t}{\epsilon},t) \, d\mathbf{x} \, dk \, d\alpha \, dt \\
&\qquad + \int_{0}^{T}\int_{0}^{2\pi}\int_{\Gamma} \sqrt{2k}\,E_{\epsilon,2}(\mathbf{x},k,\alpha,t)\,\sin\alpha\, (\D_{\alpha}\psi)(\mathbf{x},k,\alpha+\cfrac{t}{\epsilon},t) \, d\mathbf{x} \, dk \, d\alpha \, dt \\
&\quad = \int_{0}^{T}\int_{0}^{2\pi}\int_{\Gamma} \tilde{E}_{\epsilon,1}(\tilde{\mathbf{x}},t)\, \sqrt{2k}\,\cos\alpha \\
&\quad \qquad \qquad \times (\D_{\alpha}\psi)(\tilde{x}_{1}+\sqrt{2k}\,\sin\alpha, \tilde{x}_{2}-\sqrt{2k}\,\cos\alpha,k,\alpha+\cfrac{t}{\epsilon},t) \, d\tilde{\mathbf{x}} \, dk \, d\alpha \, dt  \\
&\quad \qquad + \int_{0}^{T}\int_{0}^{2\pi}\int_{\Gamma} \tilde{E}_{\epsilon,2}(\tilde{\mathbf{x}},t)\, \sqrt{2k}\,\sin\alpha \\
&\qquad \qquad \times (\D_{\alpha}\psi)(\tilde{x}_{1}+\sqrt{2k}\,\sin\alpha, \tilde{x}_{2}-\sqrt{2k}\,\cos\alpha,k,\alpha+\cfrac{t}{\epsilon},t) \, d\tilde{\mathbf{x}} \, dk \, d\alpha \, dt \, ,
\end{split}
\end{equation}
so we have
\begin{equation}
\begin{split}
&\lim_{\epsilon\,\to\,0} \int_{0}^{T}\int_{0}^{2\pi}\int_{\Gamma} \D_{k,\alpha}^{2}\phi_{\epsilon}(\mathbf{x},k,\alpha-\cfrac{t}{\epsilon},t) \, \psi(\mathbf{x},k,\alpha,t) \, d\mathbf{x} \, dk \, d\alpha \, dt \\
&\quad = \int_{0}^{T}\int_{0}^{2\pi}\int_{\Gamma} \tilde{E}_{2}(\tilde{\mathbf{x}},t)\, \cfrac{\cos\alpha}{\sqrt{2k}} \\
&\quad \qquad \qquad \times \int_{0}^{2\pi}(\D_{\alpha}\psi)(\tilde{x}_{1}+\sqrt{2k}\,\sin\alpha, \tilde{x}_{2}-\sqrt{2k}\,\cos\alpha,k,\alpha+\tau,t) \, d\tau\, d\tilde{\mathbf{x}} \, dk \, d\alpha \, dt  \\
&\quad \qquad - \int_{0}^{T}\int_{0}^{2\pi}\int_{\Gamma} \tilde{E}_{1}(\tilde{\mathbf{x}},t)\, \cfrac{\sin\alpha}{\sqrt{2k}} \\
&\quad \qquad \qquad \times \int_{0}^{2\pi} (\D_{\alpha}\psi)(\tilde{x}_{1}+\sqrt{2k}\,\sin\alpha, \tilde{x}_{2}-\sqrt{2k}\,\cos\alpha,k,\alpha+\tau,t) \, d\tau\, d\tilde{\mathbf{x}} \, dk \, d\alpha \, dt \\
&\quad = \int_{0}^{T}\int_{0}^{2\pi}\int_{\Gamma} \tilde{E}_{2}(\tilde{\mathbf{x}},t)\, \cfrac{\cos\alpha}{\sqrt{2k}} \\
&\quad \qquad \qquad \times \int_{0}^{2\pi}(\D_{\tau}\psi)(\tilde{x}_{1}+\sqrt{2k}\,\sin\alpha, \tilde{x}_{2}-\sqrt{2k}\,\cos\alpha,k,\alpha+\tau,t) \, d\tau\, d\tilde{\mathbf{x}} \, dk \, d\alpha \, dt  \\
&\quad \qquad - \int_{0}^{T}\int_{0}^{2\pi}\int_{\Gamma} \tilde{E}_{1}(\tilde{\mathbf{x}},t)\, \cfrac{\sin\alpha}{\sqrt{2k}} \\
&\quad \qquad \qquad \times \int_{0}^{2\pi} (\D_{\tau}\psi)(\tilde{x}_{1}+\sqrt{2k}\,\sin\alpha, \tilde{x}_{2}-\sqrt{2k}\,\cos\alpha,k,\alpha+\tau,t) \, d\tau\, d\tilde{\mathbf{x}} \, dk \, d\alpha \, dt \\
&\quad = 0 \, ,
\end{split}
\end{equation}
and
\begin{equation}
\begin{split}
&\lim_{\epsilon\,\to\,0} \int_{0}^{T}\int_{0}^{2\pi}\int_{\Gamma} \D_{k,\alpha}^{2}\phi_{\epsilon}(\mathbf{x},k,\alpha-\cfrac{t}{\epsilon},t) \, \psi(\mathbf{x},k,\alpha,t) \, d\mathbf{x} \, dk \, d\alpha \, dt \\
&\quad = \int_{0}^{T}\int_{0}^{2\pi}\int_{\Gamma} \tilde{E}_{1}(\tilde{\mathbf{x}},t)\, \sqrt{2k}\,\cos\alpha \\
&\quad \qquad \qquad \times \int_{0}^{2\pi}(\D_{\alpha}\psi)(\tilde{x}_{1}+\sqrt{2k}\,\sin\alpha, \tilde{x}_{2}-\sqrt{2k}\,\cos\alpha,k,\alpha+\tau,t) \, d\tau\, d\tilde{\mathbf{x}} \, dk \, d\alpha \, dt  \\
&\quad \qquad + \int_{0}^{T}\int_{0}^{2\pi}\int_{\Gamma} \tilde{E}_{2}(\tilde{\mathbf{x}},t)\, \sqrt{2k}\,\sin\alpha \\
&\quad \qquad \qquad \times \int_{0}^{2\pi} (\D_{\alpha}\psi)(\tilde{x}_{1}+\sqrt{2k}\,\sin\alpha, \tilde{x}_{2}-\sqrt{2k}\,\cos\alpha,k,\alpha+\tau,t) \, d\tau\, d\tilde{\mathbf{x}} \, dk \, d\alpha \, dt \\
&\quad = \int_{0}^{T}\int_{0}^{2\pi}\int_{\Gamma} \tilde{E}_{1}(\tilde{\mathbf{x}},t)\, \sqrt{2k}\,\cos\alpha \\
&\quad \qquad \qquad \times \int_{0}^{2\pi}(\D_{\tau}\psi)(\tilde{x}_{1}+\sqrt{2k}\,\sin\alpha, \tilde{x}_{2}-\sqrt{2k}\,\cos\alpha,k,\alpha+\tau,t) \, d\tau\, d\tilde{\mathbf{x}} \, dk \, d\alpha \, dt  \\
&\quad \qquad + \int_{0}^{T}\int_{0}^{2\pi}\int_{\Gamma} \tilde{E}_{2}(\tilde{\mathbf{x}},t)\, \sqrt{2k}\,\sin\alpha \\
&\quad \qquad \qquad \times \int_{0}^{2\pi} (\D_{\tau}\psi)(\tilde{x}_{1}+\sqrt{2k}\,\sin\alpha, \tilde{x}_{2}-\sqrt{2k}\,\cos\alpha,k,\alpha+\tau,t) \, d\tau\, d\tilde{\mathbf{x}} \, dk \, d\alpha \, dt \\
&\quad = 0 \, .
\end{split}
\end{equation}
Then we deduce that 
\begin{equation}
\nabla_{(k,\alpha)}\,\D_{\alpha} \Bigg(\int_{0}^{2\pi} \Phi(\mathbf{x},k,\alpha-\tau,\tau,t) \,d\tau \Bigg) = 0 \, ,
\end{equation}
and that the function $\displaystyle (\mathbf{x},k,\alpha) \mapsto \D_{\alpha} \Bigg(\int_{0}^{2\pi} \Phi(\mathbf{x},k,\alpha-\tau,\tau,t) \,d\tau \Bigg)$ is constant on $\Omega$. Then, for any $(\mathbf{x},k,t)$, the function $\displaystyle \alpha \mapsto \int_{0}^{2\pi} \Phi(\mathbf{x},k,\alpha-\tau,\tau,t) \,d\tau$ is monotonic and periodic on $[0,2\pi]$ and then must be constant. As a consequence, we finally obtain
\begin{equation}
\D_{\alpha} \Bigg(\int_{0}^{2\pi} \Phi(\mathbf{x},k,\alpha-\tau,\tau,t) \,d\tau \Bigg) = 0 \, .
\end{equation}
\end{proof}

~

Consequently, if we consider the hypothesis (\ref{hypothesis_E}) alongwith those from Theorem \ref{CV_GC}, we can write
\begin{equation}
\int_{0}^{2\pi} \Phi(\mathbf{x},k,\alpha-\tau,\tau,t) \,d\tau = \int_{0}^{2\pi} \Phi(\mathbf{x},k,-\tau,\tau,t) \,d\tau \, .
\end{equation}
Then, the two-scale limit model (\ref{GC_two-scale}) can be reduced to
\begin{equation} \label{GC_two-scale_reduced}
\left\{
\begin{array}{l}
\displaystyle \D_{t}G + \langle \mathcal{E}_{2} \rangle \, \D_{x_{1}}G - \langle \mathcal{E}_{1} \rangle \, \D_{x_{2}}G + \langle \mathcal{F}_{\alpha} \rangle \, \D_{\alpha} G = 0 \, , \\ \\
G(\mathbf{x},k,\alpha,0) = \cfrac{1}{2\pi} \, \tilde{f}^{0}\big(x_{1}-\sqrt{2k}\,\sin\alpha, x_{2}+\sqrt{2k}\,\cos\alpha,\sqrt{2k}\,\cos\alpha, \sqrt{2k}\,\sin\alpha\big) \, , \\ \\
\displaystyle -\nabla_{\tilde{\mathbf{x}}}\tilde{\Phi}(\tilde{\mathbf{x}},\tau,t) = \tilde{\mathcal{E}}(\tilde{\mathbf{x}},\tau,t) \, , \\ \\
\displaystyle -\Delta_{\tilde{\mathbf{x}}}\tilde{\Phi}(\tilde{\mathbf{x}},\tau,t) = \int_{\mathcal{S}} G\big(\tilde{x}_{1}+\sqrt{2k}\,\sin\alpha,\tilde{x}_{2}-\sqrt{2k}\,\cos\alpha,k,\alpha+\tau,t\big) \, dk\,d\alpha \\
\qquad \qquad \qquad \qquad \qquad \qquad \qquad \qquad \qquad \qquad \qquad \qquad \qquad \qquad - \cfrac{1}{2\pi}\,\tilde{n}_{e}(\tilde{\mathbf{x}}) \, , \\ \\
\displaystyle \mathcal{F}_{\alpha}(\mathbf{x},k,\tau,t) = \frac{\mathcal{E}_{2}(\mathbf{x},k,\alpha,\tau,t) \, \cos\alpha - \mathcal{E}_{1}(\mathbf{x},k,\alpha,\tau,t) \, \sin\alpha}{\sqrt{2k}} \, ,
\end{array}
\right.
\end{equation}
and the notation $\langle \cdot \rangle$ stands for
\begin{equation} \label{bracket_def}
\langle u \rangle(\mathbf{x},k,t) = \int_{0}^{2\pi} u(\mathbf{x},k,-\tau,\tau,t) \, d\tau \, .
\end{equation}

\begin{theorem} \label{CV-weakstar}
Under the hypotheses of Theorem \ref{CV_GC} and (\ref{hypothesis_E}), there exists a function $g$ defined on $\Gamma \times [0,T]$ such that
\begin{equation} \label{link_fg}
f(\mathbf{x},k,\alpha,t) = \frac{1}{2\pi}\,g(\mathbf{x},k,t)\,, \qquad \forall\,(\mathbf{x},k,\alpha,t) \in \Omega \times[0,T] \, ,
\end{equation}
and verifying
\begin{equation} \label{GC_weak}
\left\{
\begin{array}{l}
\displaystyle \D_{t}g + \langle \mathcal{E}_{2} \rangle \, \D_{x_{1}} g - \langle \mathcal{E}_{1} \rangle \, \D_{x_{2}}g = 0 \, ,\\ \\
\displaystyle g(\mathbf{x},k, 0) = \displaystyle \int_{0}^{2\pi} \tilde{f}^{0}(x_{1}-\sqrt{2k}\,\sin\alpha, x_{2} + \sqrt{2k}\,\cos\alpha, \sqrt{2k}\,\cos\alpha, \sqrt{2k}\,\sin\alpha) \,d\alpha \, ,\\ \\
\displaystyle \langle \mathcal{E} \rangle (\mathbf{x},k,t) = \frac{1}{2\pi} \int_{0}^{2\pi} \tilde{\mathbf{E}}\big(x_{1}-\sqrt{2k}\,\sin\alpha,x_{2}+\sqrt{2k}\,\cos\alpha,t) \, d\alpha \, , \\ \\
\displaystyle -\nabla_{\tilde{\mathbf{x}}} \tilde{\phi}(\tilde{\mathbf{x}},t) = \tilde{\mathbf{E}}(\tilde{\mathbf{x}},t) \, , \\ \\
\displaystyle -\Delta_{\tilde{\mathbf{x}}}\tilde{\phi}(\tilde{\mathbf{x}},t) = \frac{1}{2\pi} \int_{\mathcal{S}} g\big(\tilde{x}_{1}+\sqrt{2k}\,\sin\alpha,\tilde{x}_{2}-\sqrt{2k}\,\cos\alpha,k,t) \, dk\,d\alpha - \tilde{n}_{e}(\tilde{\mathbf{x}}) \, .
\end{array}
\right.
\end{equation}
\end{theorem}

\begin{proof}[Proof of Theorem \ref{CV-weakstar}]
Since $F$ and $G$ are linked by the relation (\ref{relation_FG}), $f$ does not depend on $\alpha$. Indeed, we have
\begin{equation}
\begin{split}
\D_{\alpha}f(\mathbf{x},k,\alpha,t) &= \D_{\alpha}\int_{0}^{2\pi} F(\mathbf{x},k,\alpha,\tau,t)\,d\tau \\
&= \D_{\alpha}\int_{0}^{2\pi} G(\mathbf{x},k,\alpha+\tau,t)\,d\tau \\
&= \D_{\alpha}\int_{0}^{2\pi} G(\mathbf{x},k,\tau,t) \, d\tau \\
&= 0 \, .
\end{split}
\end{equation}
Then, if we integrate (\ref{GC_two-scale_reduced}.a) in $\alpha$, we obtain
\begin{equation}
\D_{t}f + \langle \mathcal{E}_{2} \rangle \, \D_{x_{1}}f - \langle \mathcal{E}_{1} \rangle \, \D_{x_{2}}f = 0 \, ,
\end{equation}
and the equation (\ref{GC_weak}.a). \\
\indent By integrating the initial condition (\ref{GC_two-scale_reduced}.b) in $\alpha$ and dividing it by $2\pi$, we obtain the initial condition (\ref{GC_weak}.b). \\
\indent Since $\langle \mathcal{E} \rangle$ does not depend on $\alpha$, we can integrate it in $\alpha$ and divide it by $2\pi$: then, we have
\begin{equation}
\begin{split}
\langle \mathcal{E} \rangle (\mathbf{x},k,t) &= \cfrac{1}{2\pi} \, \int_{0}^{2\pi} \int_{0}^{2\pi} \mathcal{E}(\mathbf{x},k,\alpha-\tau,\tau,t) \, d\tau \, d\alpha \\
&= \cfrac{1}{2\pi} \, \int_{0}^{2\pi} \int_{0}^{2\pi} \mathcal{E}(\mathbf{x},k,\alpha,\tau,t) \, d\tau \, d\alpha \\
&= \cfrac{1}{2\pi} \, \int_{0}^{2\pi} \int_{0}^{2\pi} \tilde{\mathbf{E}} (x_{1}-\sqrt{2k}\,\sin\alpha,x_{2}+\sqrt{2k}\,\cos\alpha, t) \, d\alpha \, .
\end{split}
\end{equation}
By using similar techniques, we integrate (\ref{GC_two-scale_reduced}.d) in $\tau$ and we obtain
\begin{equation}
-\Delta_{\tilde{\mathbf{x}}}\tilde{\phi}(\tilde{\mathbf{x}},t) = \frac{1}{2\pi} \int_{\mathcal{S}} g\big(\tilde{x}_{1}+\sqrt{2k}\,\sin\alpha,\tilde{x}_{2}-\sqrt{2k}\,\cos\alpha,k,t) \, dk\,d\alpha - \tilde{n}_{e}(\tilde{\mathbf{x}}) \, .
\end{equation}
Finally, we integrate (\ref{GC_two-scale_reduced}.c) in $\tau$ and we obtain (\ref{GC_weak}.d), which concludes the proof.
\end{proof}

~

\indent We can remark that, even if they are based on the non-physical hypothesis (\ref{hypothesis_E}) for $(\tilde{\mathbf{E}}_{\epsilon})_{\epsilon\,>\,0}$, Bostan's results on the mathematical justification of the 2D finite Larmor radius approximation have been improved: indeed, we have proved that the hypotheses of Theorem \ref{CV_GC} coupled with (\ref{hypothesis_E}) are sufficient to prove Theorem \ref{CV-weakstar}. Furthermore, this result has been generalized to a non-periodic initial distribution $\tilde{f}_{0}$ and to a non-uniform electron density $\tilde{n}_{e}$.

\section{Conclusions and perspectives}

In a first part, we recalled a two-scale convergence result for a 2D Vlasov-Poisson model for a charged particle beam which is due to Fr\'enod and Sonnendr\"ucker. After adapting it to plasma modeling by adding an electron density $\tilde{n}_{e}$ and some reasonable compatibility conditions, we recalled that such a result trivially implies, up to some subsequences, the weak-* convergence of $(\tilde{f}_{\epsilon},\tilde{\mathbf{E}}_{\epsilon})$ to a couple $(\tilde{f},\tilde{\mathbf{E}})$. Then we compared this conclusion to Bostan's results which was established in \cite{Bostan_2007} under stronger assumptions. \\
\indent In a second part, we introduced a new set of variables involving the guiding-center position coordinates and the transverse part of the kinetic energy which is denoted with $k$. After rewriting the Vlasov equation in these new coordinates, we established a two-scale convergence result by using Frenod \& Sonnendr\"ucker's assumptions and the compatibility conditions for $\tilde{n}_{e}$ which were added in Theorem \ref{CV_cart}. Then, we proved that $(f_{\epsilon},\tilde{\mathbf{E}}_{\epsilon})$ weakly-* converges to a couple $(f,\tilde{\mathbf{E}})$ and, under an additional non-physical property for $(\tilde{\mathbf{E}}_{\epsilon})_{\epsilon\,>\,0}$, we established a system of constraints satisfied by $(f,\tilde{\mathbf{E}})$ through a few computations lines, remarking that this limit model is exactly the 2D finite Larmor radius model which was studied in \cite{Bostan_2007}. \\

\indent The first main remark we can do about the present work concerns the formulation of the two-scale limit model in canonical gyrokinetic coordinates under hypothesis (\ref{hypothesis_E}): indeed, the transport equations (\ref{GC_two-scale_reduced}.a) and (\ref{GC_weak}.a) indicate that $k$ is an adiabatic invariant not only for the weak-* limit system but also for the two-scale limit system. Since $k$ is the dimensionless transverse kinetic energy, this remark confirms previous results about the adiabatic invariant property of the magnetic momentum (see Littlejohn \cite{Littlejohn}, Lee \cite{Lee, Lee_2}, Brizard \textit{et al.} \cite{Brizard_PhD, Hahm-Brizard}, Grandgirard \textit{et al.} \cite{GYSELA4D}). \\
\indent The second remark which can be done is about the mathematical results which have been established in this paper: by linking previous two-scale and weak-* convergence under common assumptions, \textit{i.e.} the assumptions from Theorem \ref{CV_GC} coupled with (\ref{hypothesis_E}), we have improved Bostan's results by generalizing them to some cases involving a non-uniform electron density or a non-periodic initial distribution $\tilde{f}_{0}$. \\

\indent Since the canonical gyrokinetic variables allow us to simplify the formulation of the two-scale limit model and, under (\ref{hypothesis_E}), allow us to highlight the adiabatic invariant property of the magnetic momentum, they may be a useful tool for a mathematical justification of the full 3D finite Larmor radius model. \\
\indent From a numerical point of view, the model (\ref{GC_two-scale}) can be used to build a two-scale numerical method in order to simulate the high frequency oscillations of the solution $(\tilde{f}_{\epsilon},\tilde{\mathbf{E}}_{\epsilon})$ of the Vlasov-Poisson model (\ref{VP-2D_rescaled}), such as it has been done before for low Mach number problems (see \cite{EVA01}), charged particle beams problems (see \cite{PIC-two-scale} and \cite{Mouton_2009}) or drift problems in the ocean (see \cite{Modelling_coastal}). As an example, a numerical method based on the computation of the characteristics associated with the limit transport equation will be simpler to be developed on the formulation (\ref{GC_two-scale}.a) than the formulation (\ref{H-cart}.a): indeed, since we have 
\begin{equation}
\D_{x_{1}} \langle \mathcal{E}_{2} \rangle - \D_{x_{2}} \langle \mathcal{E}_{1} \rangle = \D_{k} \langle \mathcal{F}_{k} \rangle + \D_{\alpha} \langle \mathcal{F}_{\alpha} \rangle = 0 \, ,
\end{equation}
the transport equation (\ref{GC_two-scale}.a) can support a time splitting without losing its conservation property (see \cite{Besse-Sonnen}).

~

{\bf Acknowledgments:} We would like to thank Daniel Han-Kwan for a helpful remark he did and which allowed us to correct a mistake.

\end{document}